\documentclass[times]{kristall}	        

\usepackage{epsf}
\usepackage{amsmath,amssymb}

\DeclareMathOperator{\lcm}{lcm}
\newcommand{\CC}{\mathbb{C}}
\newcommand{\KK}{\mathbb{K}}
\newcommand{\kk}{\mathbb{L}}
\newcommand{\ZZ}{\mathbb{Z}}
\newcommand{\NN}{\mathbb{N}}
\newcommand{\RR}{\mathbb{R}}
\newcommand{\FF}{\mathbb{F}}
\newcommand{\QQ}{\mathbb{Q}\hspace{0.5pt}}
\newcommand{\OO}{\mathcal{O}}
\newcommand{\oo}{{\scriptstyle\mathcal{O}}}
\newcommand{\+}{\hfill + \hfill}

\newcommand{\qed}{\hspace*{\fill}$\square$\medskip}

\newcommand{\I}{\mathrm{i}}

\newtheorem{lem}{Lemma}
\newtheorem{prop}{Proposition}
\newtheorem{thm}{Theorem}
\newtheorem{coro}{Corollary}

\setbox\Copyrightbox\hbox{\vphantom{Yy}}
\makeatletter
\renewcommand{\@evenhead}{\underline{\hbox to 
\textwidth{{\PageFont\thepage}\hfill\MarkFont\strut \EvenRunningHead}}}%
\renewcommand{\@oddhead}{\underline{\hbox to \textwidth{{\MarkFont 
\strut\OddRunningHead}\hfill\PageFont\thepage}}}%
\renewcommand{\ps@myheadings}{}

\renewcommand{\Pedigree}{%
\begin{tabular}[t]{@{} l @{\hspace{1mm}} l @{}}
\multicolumn{2}{@{} l @{}}{\AuthorsFont\AuthorsComplete\vrule height0pt 
depth 4.83mm width0pt}\\
\raisebox{0pt}[0pt][0pt]{$^{\mathrm I}$} & \InstitutionOne\\
$^{\mathrm {II}}$ &\InstitutionTwo\vrule height0pt depth 4.83mm width0pt\strut
\end{tabular}}

\makeatother

\begin{document}
       
\maketitle		

\begin{abstract}	
\Abstract
\end{abstract}		

\bigskip
\begin{small}
\noindent
{\bf Key Words}:\hspace{\parindent} 
Lattices, Coincidence Ideals, Planar Modules, Cyclotomic Fields,
Dirichlet Series, Asymptotic Properties
\end{small}

\section{Introduction}

Given a lattice $\varGamma\subset\mathbb{R}^d$, with $d\ge 2$,
it is interesting to
know its coincidence site lattices (CSLs), which originate from
intersections of $\varGamma$ with a rotated copy. In fact,
\begin{equation}
\mathrm{SOC}(\varGamma)\; := \; \left\{ R\in\mathrm{SO}(d) \mid 
\left[\varGamma:(\varGamma\cap R\varGamma)\right]<\infty\right\}
\end{equation}
is a \emph{group}, whose structure in general (i.e., for $d>2$) is not well
understood. The coincidence index
\begin{equation}
\Sigma(R)\; = \;\left[\varGamma:(\varGamma\cap R\varGamma)\right]
\end{equation}
of a rotation $R$ is defined as the number of cosets of
$\varGamma\cap R\varGamma$ in $\varGamma$ (which can be $\infty$), and
the spectrum of finite $\Sigma$-values, 
$\Sigma\big(\mathrm{SOC}(\varGamma)\big)$, is often the first quantity
considered, followed by counting all CSLs of a given index. 
Note that one can also consider the obvious extension
to general orthogonal transformations $R\in \mathrm{O}(d)$ and the 
corresponding group $\mathrm{OC}(\varGamma)$. Coincidence lattices
play an important role in the theory of grain boundaries \cite{Bollmann,G},
and small indices can be determined in suitable experiments.

The classification of these elementary or \emph{simple}\/ coincidences
is partly done, in particular for low-dimensional lattices with
irreducible symmetries, compare \cite{BP,B} and references given there, 
but also for certain generalizations to quasicrystals \cite{PBR,B}.
The situation for various related problems \cite{BG,DCG} is similar. 
Although interesting questions concerning the 
Bravais types of the possible CSLs are still open, the general theory
in dimensions $d\le 4$ is in rather good shape, see \cite{B,P1,P2}
and references given there.

More recent is the problem of optimal lattice quantizers \cite{Sloane},
and connected with it is the question for \emph{multiple}\/
coincidences. Here, one would like to classify all lattices that can
be obtained as multiple intersections of the form
\begin{equation}
\varGamma\cap R_{1}\varGamma\cap \ldots \cap R_{m}\varGamma
\end{equation}
with $R_{\ell}\in\mathrm{SO}(d)$. One defines the corresponding
index as
\begin{equation} \label{gen-index}
   \Sigma(R_1, \ldots, R_m) \; := \;
   [\varGamma : (\varGamma \cap R_{1}\varGamma\cap 
     \ldots \cap R_{m}\varGamma)],
\end{equation}
which is finite if and only if each $R_{\ell}\in\mathrm{SOC}(\varGamma)$,
due to the mutual commensurability of the lattices $\varGamma\cap R_{\ell}
\varGamma$\/ (we shall explain this in more detail below). 
Consequently, one attaches the group 
$\big(\mathrm{SOC}(\varGamma)\big)^m$ to the setting of $m$-fold
coincidences. The corresponding spectrum is its image under $\Sigma$,
while the full (or complete) coincidence spectrum of $\varGamma$ is
\begin{equation} \label{full-spec}
    \Sigma^{}_{\varGamma} \; := \;
    \bigcup_{m\ge 1} \Sigma\big(
    (\mathrm{SOC}(\varGamma))^m\big).
\end{equation}
This is an inductive limit, with 
\[
\Sigma\big((\mathrm{SOC}(\varGamma))^m\big)\;\subset\;
\Sigma\big((\mathrm{SOC}(\varGamma))^{m'}\big)\quad \mbox{for $m \le m'$},
\] 
which often stabilizes: from a certain $m$ on, the spectra 
are stable, i.e., they do not grow any more \cite{Z}. In our examples
below, this actually happens at $m=1$, so that the spectrum is
$\Sigma^{}_{\varGamma} = \Sigma \big(\mathrm{SOC}(\varGamma)\big)$.

Clearly, multiple coincidences are also relevant in crystallography,
as they are the basis for triple or multiple junctions, in obvious
generalization of twinning, compare \cite{Gertsman}.

The problem of multiple coincidences is considerably more
involved than that of the simple ones, in particular in dimensions $d\ge 3$.
However, for $d=2$, one can rather easily extend the treatment of elementary
coincidences to multiple ones, building on the powerful and well
understood connection to the algebraic theory of cyclotomic
fields and to the analytic theory of the corresponding zeta functions.
This is precisely what we shall do below, both for lattices
and dense modules of the plane. The latter step makes the results
applicable to planar quasicrystals that are used to model the
so-called T-phases, which show a periodic stacking of planar layers
with non-crystallographic symmetries \cite{Steurer}.

Before we expand on the mathematical tools, let us set the scene with
a simple example. Afterwards, we shall treat all planar modules with
$N$-fold symmetry, class number $1$ (see below for an explanation),
and rank $\phi(N)$, where $\phi$ denotes Euler's totient function (as
defined below in Eq.~\eqref{eulerphi}). To facilitate crystallographic
applications, we also provide explicit results, both in terms of
Dirichlet series generating functions and in terms of tables; a
corresponding Mathematica notebook file is avalailable at \cite{GP}.

The following exposition is based upon previous results, and on
references \cite{PBR} and \cite{BG2} in particular. In order to 
keep our presentation concise, we shall use the notation established
there and frequently refer to these references.

\section{Example: The triangular lattice}
\label{example:triangular}

Let us consider the triangular lattice, often also called the
hexagonal lattice.  Since rotations in the plane are most 
easily written as multiplication with a complex number on the unit
circle, $e^{i\varphi}$, a suitable representative of the lattice should 
be formulated accordingly. With $\xi_{3}=e^{2\pi i/3}$, one defines
\begin{equation}
\varGamma\; = \; \left\{m+n\cdot\xi_{3} \mid m,n\in\mathbb{Z}\right\}
\; = \; \mathbb{Z}[\xi_{3}]
\label{triangular}
\end{equation}
which is a triangular lattice, with minimal distance $1$ and thus a 
fundamental cell of area $\tfrac{1}{2}\sqrt{3}$. At the same time, it 
is the set of
\emph{Eisenstein integers}\/ (also called Eisenstein-Jacobi integers),
which form the ring of integers \cite[Ch.~12.2]{HW} of the cyclotomic
field $\mathbb{Q}(\xi_{3})$, a totally complex extension of the
totally real field $\mathbb{Q}$ (and hence of degree $2$ over $\QQ$).

\begin{figure}
\epsfxsize=\columnwidth\epsfbox{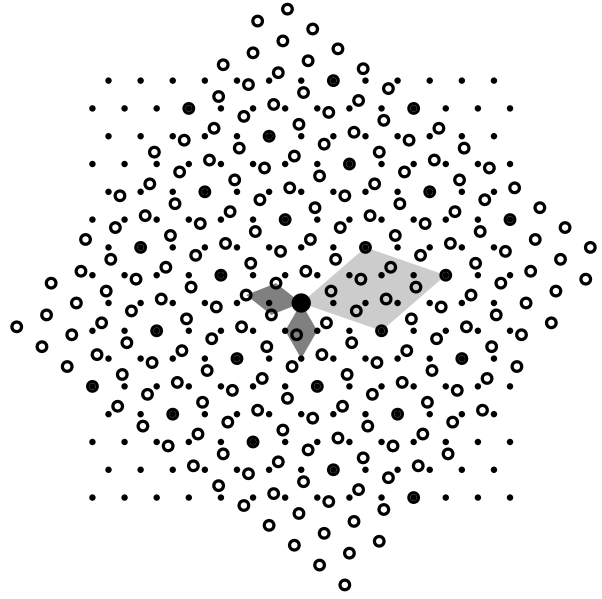}
\caption[]{Simple coincidence (with $\Sigma=7$) for the triangular lattice. 
Small dots indicate the original lattice sites, while the
circles correspond to the rotated copy. Large dots are the
coinciding sites.
Unit cells are shaded for the intersecting lattices (dark) and the 
CSL (light). 
\label{fig0}}
\end{figure}

A remarkable property of $\mathbb{Z}[\xi_{3}]$ is its unique prime
decomposition, up to units. The latter form the cyclic group
\begin{equation}
C_{6} \; = \; \left\{\pm\xi_{3}^{\ell}\mid 0\le\ell\le 2\right\}
\label{unitgroup}
\end{equation}
of rotation symmetries of $\varGamma$, with $-\xi^{}_3$ being a possible
generator for it. For any rational prime $p$ (i.e., any prime of
$\mathbb{Z}\subset\mathbb{Q}$), one of the following three possibilities
applies, compare \cite[Ch.~15.3]{HW}:\smallskip
\begin{enumerate}
\item $p=3$; this prime is called \emph{ramified}, and factorizes as 
  $3=(1-\xi_{3})(1-\bar{\xi}_{3})=-\bar{\xi}_{3}(1-\xi_{3})^2$. Up to a
  unit, $3$ is thus the square of a prime in $\mathbb{Z}[\xi_{3}]$.
\item $p\equiv 2 \; (3)$; these primes are called \emph{inert}, because 
  they stay prime in $\mathbb{Z}[\xi_{3}]$.
\item $p\equiv 1 \; (3)$; these are the (complex) \emph{splitting}\/ primes,
  because they factorize as $p=\omega_{p}\cdot\bar{\omega}_{p}$ into a pair of
  complex conjugate primes of $\mathbb{Z}[\xi_{3}]$ that are not
  associated to one another (meaning that $\omega_{p}/\bar{\omega}_{p}$
  is not a unit in $\ZZ[\xi_3]$).
\end{enumerate}\smallskip
The last type of primes is the key to solving the coincidence
problem. As was shown in \cite{PBR}, any simple coincidence
rotation of $\varGamma$ of (\ref{triangular}) is of the form
\begin{equation}
e^{i\varphi} \; = \; \varepsilon\cdot\prod_{p\equiv 1 \, (3)} 
\left(\frac{\omega_{p}}{\bar{\omega}_{p}}\right)^{t_{p}}
\label{trirot}
\end{equation}
where the product runs over all rational primes $\equiv 1 \, (3)$,
with all $t_{p}\in\mathbb{Z}$ (only \emph{finitely}\/ many of them
differing from $0$), and where $\varepsilon$ is a unit in $\ZZ[\xi_3]$,
i.e., $\varepsilon\in C_6$ from Eq.~\eqref{unitgroup}. An example with
the smallest non-trivial index is shown in Figure~\ref{fig0}.

The meaning (and the basic ingredient of the proof) is that
$e^{i\varphi}$ of \eqref{trirot} rotates one Eisenstein integer (the
denominator) into another (the numerator).  The corresponding CSL is
then generated by $w$ and $\xi_{3}\cdot w$, where
\begin{equation}
   w \; =
  \prod_{\stackrel{\scriptstyle p\equiv 1 \, (3)}{\scriptstyle t_{p}>0}} 
  (\omega_{p})^{t_{p}} \cdot
  \prod_{\stackrel{\scriptstyle p\equiv 1 \, (3)}{\scriptstyle t_{p}<0}} 
  (\bar{\omega}_{p})^{-t_{p}}
\label{omega}
\end{equation}
is the numerator of $e^{i\varphi}$ in \eqref{trirot}.

Clearly, the CSL that emerges from Eq.~\eqref{omega} is the principal
ideal $(w) = w \, \ZZ[\xi_3]$. Its index in $\ZZ[\xi_3]$, also
called its total norm, is the corresponding coincidence index and
reads
\begin{equation}
\Sigma(\varphi) \; = \prod_{p\equiv 1 \, (3)} p^{|t_{p}|}
\end{equation}
with $\Sigma=1$ precisely for $e^{i\varphi}=\varepsilon\in C_{6}$.

This shows that the group of coincidence rotations has the structure
\begin{equation}
  \mathrm{SOC}(\varGamma) \; = \;
  \mathrm{SOC}(\ZZ[\xi_3]) \; = \;
  C_{6}\times\mathbb{Z}^{(\aleph_{0})},
\label{3group}
\end{equation}
where $\mathbb{Z}^{(\aleph_{0})}$ stands for the direct sum of
countably many infinite cyclic groups, each with a generator of
the form $\omega_{p}/\bar{\omega}_{p}$, and the restriction
mentioned after Eq.~\eqref{trirot}. The set of possible indices,
$\Sigma\big(\mathrm{SOC}(\ZZ[\xi_3])\big)$, forms a semigroup with
unit, generated by $1$ and the rational primes $p\equiv 1 \, (3)$, the
latter being called \emph{basic indices}.

If $c^{}_3 (k)$ denotes the number of CSLs of $\varGamma$ of index $k$,
it is a multiplicative arithmetic function (i.e., $c^{}_3 (k\ell)=
c^{}_3 (k) c^{}_3 (\ell)$
whenever $k,\ell$ are coprime). Consequently, its
determination is most easily achieved by means of a Dirichlet
series generating function. The result reads \cite{PBR}
\begin{align}
\varPhi_{3}(s) & \; = \; \sum_{k=1}^{\infty}\frac{c^{}_3 (k)}{k^{s}}
\; =\prod_{p\equiv 1 \, (3)}\frac{1+p^{-s}}{1-p^{-s}}\\
& \; = \;
1\!+\!\tfrac{2}{7^s}\!+\!\tfrac{2}{13^s}\!+\!\tfrac{2}{19^s}\!+\!
  \tfrac{2}{31^s}\!+\!\tfrac{2}{37^s}\!+\!\tfrac{2}{43^s}\!+\!\tfrac{2}{49^s}
+\ldots\notag
\end{align}
where terms with numerators $>2$ emerge as suitable composite
denominators show up.
This can also be expressed in terms of zeta functions,
\begin{equation}
   \varPhi_{3}(s) \; = \;
   \frac{\zeta^{}_{\mathbb{Q}(\xi_{3})}(s)}
        {(1+3^{-s})\, \zeta(2s)},
\label{phi3}
\end{equation}
where $\zeta(s)=\sum_{k=1}^{\infty}k^{-s}$ is Riemann's
zeta function, and 
\begin{align}
  \zeta^{}_{\mathbb{Q}(\xi_{3})}(s) 
  & \; = \; \sum_{k=1}^{\infty} \frac{a^{}_{3} (k)}{k^s} \\
  & \; = \; \frac{1}{1-3^{-s}}
  \prod_{p\equiv 1 \, (3)} \frac{1}{(1-p^{-s})^2}
  \prod_{p\equiv 2 \, (3)} \frac{1}{1-p^{-2s}} \notag 
\end{align}
is Dedekind's zeta function of the cyclotomic field $\QQ(\xi_{3})$, 
with $a^{}_{3} (k)$ the number of ideals of $\ZZ[\xi_3]$ of index $k$,
see \cite{W} for details. Both zeta functions converge absolutely
on the right half-plane $\{ s\in\CC \mid \mathrm{Re}(s)>1\}$, simply written
as $\{ \mathrm{Re}(s)>1\}$ from now on.

Since all CSLs of $\varGamma$ are themselves triangular lattices (in
particular, they are similarity sublattices of $\varGamma$), the
problem of simple coincidences is thus solved.
{}From Eq.~\eqref{phi3}, one can also extract the asymptotic behaviour
of the arithmetic function $c^{}_{3}(k)$. The result is \cite{PBR}
\begin{equation}
\sum_{k\le x} c^{}_{3}(k) \;\sim\; 
x \cdot \big( \mbox{res}_{s=1}\,\varPhi_{3}(s)\big)
\;=\; x\cdot \frac{\sqrt{3}}{2\pi},
\label{simpleasymp3} 
\end{equation}
where $f(x)\sim g(x)$ means $\lim_{x\to\infty}\frac{f(x)}{g(x)}=1$.

Let us now consider \emph{multiple}\/ coincidences (which are meant
to include the simple ones, of course). The starting point is the 
rather obvious observation that
\begin{equation}
\label{eq:inter}
  \varGamma\cap R_{1}\varGamma\cap\ldots\cap R_{m}\varGamma \; = \;
  \bigcap_{\ell=1}^{m}(\varGamma\cap R_{\ell}\varGamma)\,.
\end{equation}
If the result is still a sublattice of $\varGamma$ (and thus of
finite index in it), each $(\varGamma\cap R_{\ell}\varGamma)$ itself
must be a CSL. This necessary condition is also sufficient, because
any two sublattices of $\varGamma$ are commensurate, i.e., their
intersection is still a (full) sublattice of $\varGamma$. So, we
get a CSL in \eqref{eq:inter} if and only if all $R_{\ell}$ are
elements of $\mathrm{SOC}(\varGamma)$.

As follows immediately from Eq.~\eqref{omega}, 
each simple CSL is an ideal of $\ZZ[\xi_3]$ of the
form $(w) := w \, \ZZ[\xi_3]$, i.e., a principal ideal.
Note that all ideals of $\ZZ[\xi_3]$ are principal, but not all
of them appear as CSLs. Any multiple CSL is an intersection of 
simple CSLs. Consequently, one has
\begin{equation}
\bigcap_{i=1}^{m}(w^{}_{i}) \; = \;
\big(\mathrm{lcm}\{w^{}_{1},\ldots,w^{}_{m}\}\big),
\label{lcm-eq}
\end{equation}
which is indeed a principal ideal again,
due to the unique prime factorization (up to units) in
$\ZZ[\xi_{3}]$. Since each $w_{i}$ is of the form
\eqref{omega}, one quickly checks that multiple coincidences cannot
enlarge the list of possible indices, though they \emph{can}\/ lead to
more solutions for a given index, hence to genuinely multiple CSLs.

Grouping the solutions from Eq.~\eqref{lcm-eq} according to their
coincidence indices, one finds the following
\begin{prop}\label{prop1}
   If\/ $k=p_{1}^{r_{1}}\cdot\ldots\cdot p_{t}^{r_{t}}$ is the prime
   decomposition of an integer\/ $k>1$ into rational primes, the triangular 
   lattice\/ $\varGamma$ has no CSL of index\/ $k$ unless \emph{all}\/ 
   $p_{j}\equiv 1 \; (3)$. 
   In that case, the number of multiple CSLs of index\/ $k$ is given
   by the multiplicative arithmetic function
\[
  b^{}_3 (k) \; = \; \prod_{j=1}^{t}(r_{j}+1) .
\]
Moreover, the Euler product representation of the corresponding
Dirichlet series generating function reads
\[
  \varPsi_{3}(s) \; = \; \sum_{k=1}^{\infty} \frac{b^{}_3 (k)}{k^s}
   \; = \prod_{p\equiv 1 \, (3)}\frac{1}{(1-p^{-s})^2},
\]
with absolute convergence on\/ $\{\mathrm{Re}(s)>1\}$. 
\end{prop}
\noindent {\sc Proof}.
Clearly, each $(w_{i})$ in \eqref{lcm-eq} must be a simple CSL of
$\varGamma$ for the total coincidence index to be finite, hence all
$p_{j}\equiv 1 \, (3)$ from Eq.~\eqref{trirot}. 

Then, each factor $p^r$ of $k$, as $p=\omega\bar{\omega}$ is a complex
splitting prime, can contribute a factor of the form
$(\omega^{r-\ell}\bar{\omega}^{\ell})$ to the lcm of
Eq.~\eqref{lcm-eq}, with $0\le\ell\le r$. Any principal ideal
$(\omega^{r-\ell}\bar{\omega}^{\ell})$, in turn, is either a
simple CSL itself (if $\ell=0$ or $\ell=r$) or the intersection of 
two simple CSLs.

So, for the prime $p$, which enters the decomposition of
$k$ as $p^r$, this amounts to $r+1$ different possibilities. 
Since the different primes give independent choices of this kind, 
one obtains $b^{}_{3} (k)$ as the product stated above, and the
multiplicativity of this arithmetic function is obvious.

The Dirichlet series has the form claimed, as one can see from its
Euler factors, using the identity 
$(1-x)^{-2}=\sum_{\ell\ge 0}(\ell+1)x^{\ell}$ with
$x=p^{-s}$; the convergence result is standard.
\qed\smallskip

A comparison of $\varPhi_{3}$ and $\varPsi_{3}$ shows that
\[
  \varPsi_{3}(s)\!-\!\varPhi_{3}(s) \: = \: 
  \tfrac{1}{49^s}\!+\!\tfrac{1}{169^s}\!+\!\tfrac{2}{343^s}\!+\! 
  \tfrac{1}{361^s}\!+\!\tfrac{2}{637^s}\!+\!\tfrac{2}{931^s}\!+\! 
  \ldots
\]
encapsulates the statistics of the genuinely multiple CSLs.  In
particular, a comparison with $\zeta^{}_{\mathbb{Q}(\xi_{3})}(s)$
shows that all ideals of norm $k$, with the additional condition that
$k$ factorizes into rational primes $p\equiv 1 \; (3)$ only, are CSLs
for multiple (in fact, single or double) intersections.  Consequently,
$\varPsi_{3}(s)$ is also the generating function for double
intersections alone, and one has
\begin{align}
\label{trispec}
\Sigma_{\varGamma} &\; = \; \Sigma\big(\mathrm{SOC}(\varGamma)\big)\\
&\; = \; \Big\{ \prod_{i=1}^{t} p_{i}^{r_{i}} \mid 
\mbox{all $p_{i}^{}\equiv 1\, (3)$, $r_{i}\in\NN$ and $t\ge 0$}\Big\}. \notag 
\end{align}

\begin{figure}
\epsfxsize=\columnwidth\epsfbox{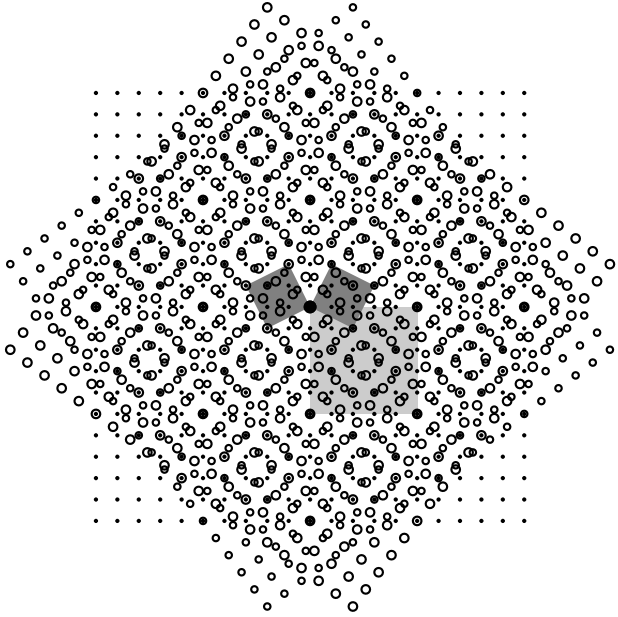}
\caption[]{Double coincidence (with $\Sigma=25$) for the square lattice. 
Unit cells are shaded for simple (dark) and double (light) CSLs. Small
dots indicate the original lattice sites, while the two types of
circles correspond to the two rotated copies. Large dots are formed by
double coincidences, i.e., at sites which belong to all three
lattices.
\label{fig1}}
\end{figure}

To expand on the first difference between simple and double
coincidences, consider $p=7$, which is the smallest rational prime
$\equiv 1 \; (3)$. It splits as $7=\omega\cdot\bar{\omega}$ with
$\omega=2-\xi_{3}$, and both the ideals $(\omega^{2})$ and
$(\bar{\omega}^2)$ appear as simple CSLs (of index $\Sigma=49$),
while $(\omega\bar{\omega})=(7)$ is only possible as multiple CSL --
namely as the intersection of the ideals $(\omega)$ and
$(\bar{\omega})$, both being simple CSLs (of index $7$, as in
Figure~\ref{fig0}). Note that $(\omega\bar{\omega}) = (\omega)\cap
(\bar{\omega})$ is, at the same time, a simple CSL for the lattice
$(\omega)$, because the latter is rotated into $(\bar{\omega})$ by
$e^{i\varphi} = \bar{\omega}/\omega$.

Asymptotically, the surplus of multiple CSLs leads 
to an additional factor in comparison to Eq.~\eqref{simpleasymp3}, so that
\begin{equation}
  \sum_{k\le x} b^{}_3 (k) \;\sim\; 
   x \cdot \big( \mbox{res}_{s=1}\,\varPsi_{3}(s)\big)
   \;=\; x\cdot \frac{\sqrt{3}}{2\pi} \cdot q^{}_{3}
\end{equation}
with 
\begin{equation}
  q^{}_{3} \; := \; \lim_{s\to 1} \frac{\varPsi_3 (s)}{\varPhi_3 (s)}
  \; = \prod_{p\equiv 1\, (3)} \frac{p^2}{p^2-1} \; \simeq \; 1.034\,015.
\end{equation}
The convergence of this product turns out to be rather slow (to put it mildly),
whence one needs a different method to actually calculate
$q^{}_{3}$. As we shall show below in greater generality, this is
achieved by means of the identity
\[
  q^{}_{3} \; = \;  \prod_{\ell=1}^{\infty} 
\big(\varPhi_3 (2^\ell)\big)^{1/2^\ell}.
\]
Although this looks just like another infinite product representation, 
its convergence is spectacularly fast\footnote{We thank Peter Pleasants for 
pointing this out to us.}. Also, it is monotonically increasing, which is 
helpful to derive error bounds for the actual numerical calculations
in this case.

\section{Cyclotomic integers with class number one}

The example of the triangular lattice was chosen because it is both
simple and paradigmatic. The square lattice, when identified with the
ring $\mathbb{Z}[i]$ of Gaussian integers \cite{HW}, can be treated in
the same way, leading to analogous expressions, with $C_{4}=\{\pm 1,
\pm i\}$ instead of $C_{6}$ and with the congruence condition $p\equiv
1 \; (4)$ rather than $p\equiv 1 \; (3)$, see \cite{B} for
details. For the square lattice, the first example of a genuinely
multiple CSL occurs for $\Sigma=25$. The corresponding situation is
illustrated in Figure~\ref{fig1}, and is completely analoguous to the
situation met above for the triangular lattice, with $5=(1+2i)(1-2i)$.

As we shall now show, the situation actually extends immediately to
all rings of integers $\mathbb{Z}[\xi_{n}]$ in cyclotomic fields
$\mathbb{Q}(\xi_{n})$ with class number $1$ (excluding $\QQ$ itself).
These emerge for the following 29 choices of $n$,
\begin{equation}
\begin{array}{lcl}
n&\in&\{3,4,5,7,8,9,11,12,13,15,16,\\
  &&\hphantom{\{}17,19,20,21,24,25,27,28,32,\\
  &&\hphantom{\{}33,35,36,40,44,45,48,60,84\}.
\end{array}
\label{liste}
\end{equation}
Here, we view each set $\OO_{n}=\mathbb{Z}[\xi_{n}]$ (i.e.,
the ring of polynomials in $\xi_n$) as a point set in $\CC$, where
$\xi_{n}$ is a primitive $n$-th root of unity. To be explicit (which is
not necessary), we choose $\xi_{n}=\exp(2\pi\I/n)$. Apart from $n=1$
(where $\OO_{1}=\mathbb{Z}$ is one-dimensional), the values
of $n$ in \eqref{liste} correspond to all cases where
$\mathbb{Z}[\xi_{n}]$ is a principal ideal domain and thus has class
number one, see \cite{W,BG} for details. If $n$ is odd, we have
$\OO_{n} = \OO_{2n}$.  Consequently,
$\OO_{n}$ has $N$-fold symmetry, where
\begin{equation}
N\; =\; N(n)\;=\; \lcm (2,n) .
\label{symm}
\end{equation}
To avoid duplication of results, integers $n\equiv 2\; (4)$ do neither
appear in the above list \eqref{liste} nor in our further exposition.

The values of $n$ from the list \eqref{liste} are naturally grouped 
according to $\phi(n)$, which is Euler's totient function
\begin{equation} \label{eulerphi}
   \phi(n) \; = \; \big\lvert 
   \{ 1\le k \le n \mid \gcd(k,n)=1 \} \big\rvert\,.
\end{equation}
The set $\OO_{n}$, which is the so-called maximal  order \cite{BS} 
of the cyclotomic field $\QQ(\xi_n)$, turns into a lattice in
$\RR^{\phi(n)}$ by means of Minkowski's embedding, compare
\cite{BS}, while it is a dense subset of $\RR^2$ whenever
$\phi(n)>2$.
Note that $\phi(n)=2$ covers the two crystallographic cases $n=3$
(triangular lattice) and $n=4$ (square lattice), while
\mbox{$\phi(n)=4$} means $n\in\{5,8,12\}$ which are the standard
symmetries of planar quasicrystals, as they occur as layers in
so-called T-phases, compare \cite{Steurer}. Again, $n=6$ and $n=10$ 
are covered implicitly, as explained above.

The rings of integers $\OO_{n}$ arise in quasicrystal theory in
several ways. One is in the form of the Fourier module that supports the
Bragg peaks of X-ray diffraction \cite{Steurer}.  Another is via the
limit translation module of a (discrete) quasiperiodic tiling
\cite{B2}. Characteristic points of the latter (e.g., vertex points)
are then model sets (or cut and project sets) on the basis of the
entire ring $\OO_n$, viewed in parallel as a (dense) point set in
$\RR^2$ and as a lattice in $\RR^{\phi(n)}$, see \cite{B,P} for
details. Such sets are also called cyclotomic model sets. 

To come to an efficient formulation, one needs another field, namely
the maximal real subfield of $\QQ(\xi_{n})$, and its ring of integers.
{}From now on, we shall thus use the following notation:
\begin{equation}
\begin{array}{rclcrcl}
 \KK_{n} & \!=\! & \QQ(\xi_{n}),  & &
 \kk_{n} & \!=\! & \QQ(\xi_{n}+\bar{\xi}_{n}), \\
 \OO_{n} & \!=\! & \ZZ[\xi_{n}],  & &
 \oo_{n} & \!=\! & \ZZ[\xi_{n}+\bar{\xi}_{n}]. 
\end{array}
\end{equation}
The next result about the relation between these fields and rings
of integers is standard \cite{W}.

\begin{lem}
  Let $\xi_{n}$ be a primitive $n$-th root of unity. The smallest
  field extension of\/ $\QQ$ that contains $\xi_{n}$ is the cyclotomic
  field\/ $\KK_{n}$. It is an extension of\/ $\QQ$ of degree
  $\phi(n)$, with\/ $\OO_{n}$ as its ring of integers\/ $($and
  its maximal order\/$)$.

Moreover, $\kk_{n}$ is the maximal real subfield of\/ $\KK_{n}$, with\/
$\oo_{n}$ as its ring of integers. For $n\ge 3$, $\KK_{n}$ is an
extension of\/ $\kk_{n}$ of degree\/ $2$, while\/ $\kk_{n}$ has degree\/
$\tfrac{1}{2}\phi(n)$ over\/ $\QQ$.\qed
\end{lem}

In our example in Section~\ref{example:triangular}, we had $n=3$.
Consequently, with $\xi_{3}=e^{2\pi i/3}$, one finds
$\xi_{3}+\bar{\xi}_{3}=-1$, so that $\kk_{3}=\QQ$ and $\oo_{3}=\ZZ$.
Also, $\kk_4=\QQ$ and $\oo_4=\ZZ$, while for $\phi(n)=4$, one finds
the three cases
\[
\begin{array}{ccccccc}
 \kk_5 & \!=\! & \QQ(\sqrt{5}) ,  && \oo_5 & \!=\! & \ZZ[\tau], \\
 \kk_8 & \!=\! & \QQ(\sqrt{2}) ,  && \oo_8 & \!=\! & \ZZ[\sqrt{2}],\\
 \kk_{12} & \!=\! & \QQ(\sqrt{3}) , && \oo_{12} & \!=\! & \ZZ[\sqrt{3}],
\end{array}
\]
with $\tau=(\sqrt{5}+1)/2$ being the golden ratio.

\section{Simple coincidences for cyclotomic integers}

One advantage of the number theoretic formulation is that the
treatment of $\varGamma=\OO_{3}$ from Section~\ref{example:triangular}
can be extended to $\OO_{n}$ for all $n$ from the list \eqref{liste}.
For simple coincidences, this was done by Pleasants et al.\ in
\cite{PBR}. This first requires the extension of the basic definitions
from lattices to more general $\ZZ$-modules of finite rank, imbedded
in Euclidean space. This is possible without difficulty in the natural
group theoretic setting \cite{B,PBR}, leading to the concept of a
coincidence site module (CSM).

\begin{table*}
\caption[]{Complex splitting primes of
$\KK_n/\QQ$, with $n$ from the list \eqref{liste}, that originate from 
rational primes $p\nmid n$ (unramified). 
The symbol $p^{\ell_p}_{(k)}$ means that primes 
$p\equiv k\; (n)$ contribute via $p^{\ell_p}$ as basic 
index, where $\ell_p$ is the smallest integer\vphantom{$k^\ell_(m)$} such that
$k^{\ell_p}\equiv 1\; (n)$, and the integer $m_p$ that appears
in Theorems \ref{thm1} and \ref{thm2} is\vphantom{$k^\ell_m$}
$m_p=\phi(n)/\ell_p$. Note that $m_p$ is always even.
\label{tab:basic}}
\def\Strut{\large\strut}
\renewcommand{\arraystretch}{1.5}
\begin{tabular*}{\textwidth}{@{\extracolsep{\fill}} @{}crl@{}}  
 \hline
\vspace{-0.5pt}$\phi(n)$ & $n$ & basic indices \Strut\\
 \hline
 2 &  3 &  $p_{(1)}^{1}$ \Strut\\
   &  4 &  $p_{(1)}^{1}$ \rule[-2ex]{0ex}{2ex}\Strut\\
 \hline
 4 &  5 &  $p_{(1)}^{1}$ \Strut\\
   &  8 &  $p_{(1)}^{1}$, $p_{(3)}^{2}$, $p_{(5)}^{2}$ \Strut\\
   & 12 &  $p_{(1)}^{1}$, $p_{(5)}^{2}$, $p_{(7)}^{2}$ 
\rule[-2ex]{0ex}{2ex}\Strut\\
 \hline
 6 &  7 &  $p_{(1)}^{1}$, $p_{(2)}^{3}$, $p_{(4)}^{3}$ \Strut\\
   &  9 &  $p_{(1)}^{1}$, $p_{(4)}^{3}$, $p_{(7)}^{3}$
\rule[-2ex]{0ex}{2ex}\Strut\\
 \hline 
 8 & 15 &  $p_{(1)}^{1}$, $p_{(4)}^{2}$, $p_{(11)}^{2}$, $p_{(2)}^{4}$, 
           $p_{(7)}^{4}$, $p_{(8)}^{4}$, $p_{(13)}^{4}$ \Strut\\
   & 16 &  $p_{(1)}^{1}$, $p_{(7)}^{2}$, $p_{(9)}^{2}$, $p_{(3)}^{4}$, 
           $p_{(5)}^{4}$, $p_{(11)}^{4}$, $p_{(13)}^{4}$ \Strut\\
   & 20 &  $p_{(1)}^{1}$, $p_{(9)}^{2}$, $p_{(11)}^{2}$, $p_{(3)}^{4}$, 
           $p_{(7)}^{4}$, $p_{(13)}^{4}$, $p_{(17)}^{4}$ \Strut\\
   & 24 &  $p_{(1)}^{1}$, $p_{(5)}^{2}$, $p_{(7)}^{2}$, $p_{(11)}^{2}$, 
           $p_{(13)}^{2}$, $p_{(17)}^{2}$, $p_{(19)}^{2}$ 
\rule[-2ex]{0ex}{2ex}\Strut\\
\hline
10 & 11 &  $p_{(1)}^{1}$, $p_{(3)}^{5}$, $p_{(4)}^{5}$, $p_{(5)}^{5}$, 
           $p_{(9)}^{5}$ 
\rule[-2ex]{0ex}{2ex}\Strut\\ 
\hline
12 & 13 &  $p_{(1)}^{1}$, $p_{(3)}^{3}$, $p_{(9)}^{3}$ \Strut\\
   & 21 &  $p_{(1)}^{1}$, $p_{(8)}^{2}$, $p_{(13)}^{2}$, $p_{(4)}^{3}$,
           $p_{(16)}^{3}$, $p_{(2)}^{6}$, $p_{(10)}^{6}$, $p_{(11)}^{6}$, 
           $p_{(19)}^{6}$ \Strut\\
   & 28 &  $p_{(1)}^{1}$, $p_{(13)}^{2}$, $p_{(15)}^{2}$, $p_{(9)}^{3}$, 
           $p_{(25)}^{3}$, $p_{(5)}^{6}$, $p_{(11)}^{6}$, $p_{(17)}^{6}$, 
           $p_{(23)}^{6}$ \Strut\\
   & 36 &  $p_{(1)}^{1}$, $p_{(17)}^{2}$, $p_{(19)}^{2}$, $p_{(13)}^{3}$, 
           $p_{(25)}^{3}$, $p_{(5)}^{6}$, $p_{(7)}^{6}$, $p_{(29)}^{6}$, 
           $p_{(31)}^{6}$
\rule[-2ex]{0ex}{2ex}\Strut\\
\hline
16 & 17 &  $p_{(1)}^{1}$ \Strut\\
   & 32 &  $p_{(1)}^{1}$, $p_{(15)}^{2}$, $p_{(17)}^{2}$, $p_{(7)}^{4}$, 
           $p_{(9)}^{4}$, $p_{(23)}^{4}$, $p_{(25)}^{4}$, $p_{(3)}^{8}$, 
           $p_{(5)}^{8}$, $p_{(11)}^{8}$, $p_{(13)}^{8}$, $p_{(19)}^{8}$, 
           $p_{(21)}^{8}$, $p_{(27)}^{8}$, $p_{(29)}^{8}$ \Strut\\
   & 40 &  $p_{(1)}^{1}$, $p_{(9)}^{2}$, $p_{(11)}^{2}$, $p_{(19)}^{2}$,
           $p_{(21)}^{2}$, $p_{(29)}^{2}$, $p_{(31)}^{2}$, $p_{(3)}^{4}$, 
           $p_{(7)}^{4}$, $p_{(13)}^{4}$, $p_{(17)}^{4}$, $p_{(23)}^{4}$,
           $p_{(27)}^{4}$, $p_{(33)}^{4}$, $p_{(37)}^{4}$ \Strut\\
   & 48 &  $p_{(1)}^{1}$, $p_{(7)}^{2}$, $p_{(17)}^{2}$, $p_{(23)}^{2}$,
           $p_{(25)}^{2}$, $p_{(31)}^{2}$, $p_{(41)}^{2}$, $p_{(5)}^{4}$,
           $p_{(11)}^{4}$, $p_{(13)}^{4}$, $p_{(19)}^{4}$, $p_{(29)}^{4}$,
           $p_{(35)}^{4}$, $p_{(37)}^{4}$, $p_{(43)}^{4}$ \Strut\\
   & 60 &  $p_{(1)}^{1}$, $p_{(11)}^{2}$, $p_{(19)}^{2}$, $p_{(29)}^{2}$, 
           $p_{(31)}^{2}$, $p_{(41)}^{2}$, $p_{(49)}^{2}$, $p_{(7)}^{4}$, 
           $p_{(13)}^{4}$, $p_{(17)}^{4}$, $p_{(23)}^{4}$, $p_{(37)}^{4}$, 
           $p_{(43)}^{4}$, $p_{(47)}^{4}$, $p_{(53)}^{4}$
\rule[-2ex]{0ex}{2ex}\Strut\\
\hline
18 & 19 &  $p_{(1)}^{1}$, $p_{(7)}^{3}$, $p_{(11)}^{3}$, $p_{(4)}^{9}$, 
           $p_{(5)}^{9}$, $p_{(6)}^{9}$, $p_{(9)}^{9}$, $p_{(16)}^{9}$, 
           $p_{(17)}^{9}$ \Strut\\
   & 27 &  $p_{(1)}^{1}$, $p_{(10)}^{3}$, $p_{(19)}^{3}$, $p_{(4)}^{9}$, 
           $p_{(7)}^{9}$, $p_{(13)}^{9}$, $p_{(16)}^{9}$, $p_{(22)}^{9}$, 
           $p_{(25)}^{9}$
\rule[-2ex]{0ex}{2ex}\Strut\\ 
\hline
20 & 25 &  $p_{(1)}^{1}$, $p_{(6)}^{5}$, $p_{(11)}^{5}$, $p_{(16)}^{5}$, 
           $p_{(21)}^{5}$ \Strut\\
   & 33 &  $p_{(1)}^{1}$, $p_{(10)}^{2}$, $p_{(23)}^{2}$, $p_{(4)}^{5}$, 
           $p_{(16)}^{5}$, $p_{(25)}^{5}$, $p_{(31)}^{5}$, $p_{(5)}^{10}$, 
           $p_{(7)}^{10}$, $p_{(13)}^{10}$, $p_{(14)}^{10}$, $p_{(19)}^{10}$, 
           $p_{(20)}^{10}$, $p_{(26)}^{10}$, $p_{(28)}^{10}$ \Strut\\ 
   & 44 &  $p_{(1)}^{1}$, $p_{(21)}^{2}$, $p_{(23)}^{2}$, $p_{(5)}^{5}$, 
           $p_{(9)}^{5}$, $p_{(25)}^{5}$, $p_{(37)}^{5}$, $p_{(3)}^{10}$, 
           $p_{(13)}^{10}$, $p_{(15)}^{10}$, $p_{(17)}^{10}$, $p_{(27)}^{10}$, 
           $p_{(29)}^{10}$, $p_{(31)}^{10}$, $p_{(41)}^{10}$ 
\rule[-2ex]{0ex}{2ex}\Strut\\
\hline
24 & 35 &  $p_{(1)}^{1}$, $p_{(6)}^{2}$, $p_{(29)}^{2}$, $p_{(11)}^{3}$, 
           $p_{(16)}^{3}$, $p_{(8)}^{4}$, $p_{(13)}^{4}$, $p_{(22)}^{4}$, 
           $p_{(27)}^{4}$, $p_{(4)}^{6}$, $p_{(9)}^{6}$, $p_{(26)}^{6}$, 
           $p_{(31)}^{6}$, \Strut\\
   &    &  $p_{(2)}^{12}$, $p_{(3)}^{12}$, $p_{(12)}^{12}$,
           $p_{(17)}^{12}$, $p_{(18)}^{12}$, $p_{(23)}^{12}$, $p_{(32)}^{12}$,
           $p_{(33)}^{12}$ \Strut\\
   & 45 &  $p_{(1)}^{1}$, $p_{(19)}^{2}$, $p_{(26)}^{2}$, $p_{(16)}^{3}$,
           $p_{(31)}^{3}$, $p_{(8)}^{4}$, $p_{(17)}^{4}$, $p_{(28)}^{4}$,
           $p_{(37)}^{4}$, $p_{(4)}^{6}$, $p_{(11)}^{6}$, $p_{(34)}^{6}$,
           $p_{(41)}^{6}$, \Strut\\
   &    &  $p_{(2)}^{12}$, $p_{(7)}^{12}$, $p_{(13)}^{12}$, 
           $p_{(22)}^{12}$, $p_{(23)}^{12}$, $p_{(32)}^{12}$, $p_{(38)}^{12}$,
           $p_{(43)}^{12}$ \Strut\\
   & 84 &  $p_{(1)}^{1}$, $p_{(13)}^{2}$, $p_{(29)}^{2}$, $p_{(41)}^{2}$, 
           $p_{(43)}^{2}$, $p_{(55)}^{2}$, $p_{(71)}^{2}$, $p_{(25)}^{3}$, 
           $p_{(37)}^{3}$, $p_{(5)}^{6}$, $p_{(11)}^{6}$, $p_{(17)}^{6}$,
           $p_{(19)}^{6}$, \Strut\\
   &    &  $p_{(23)}^{6}$, $p_{(31)}^{6}$, $p_{(53)}^{6}$,
           $p_{(61)}^{6}$, $p_{(65)}^{6}$, $p_{(67)}^{6}$, $p_{(73)}^{6}$,
           $p_{(79)}^{6}$
\rule[-2ex]{0ex}{2ex}\Strut\\
\hline
\end{tabular*}
\renewcommand{\arraystretch}{1}
\end{table*}

Here, we shall first summarise the main results of
\cite{PBR}, and then completely work out all $29$ cases. The extension
to multiple coincidences follows in the next section.

Below, we shall always make the assumption that we deal
with $\OO_n$, where $n$ is an element from our list \eqref{liste}.
Since these are the cases with class number one, we shall abbreviate
this assumption as {\bf (CN1)}. In this case, also the maximal real
subfields $\kk_{n}$ have class number one \cite[Prop.~11.19]{W}.

As in our example ($n=3$) in Section~\ref{example:triangular}, a prime
$\mathfrak{p}$ in $\oo_n$, relative to $\OO_n$, is either ramified,
inert, or splits into a pair of complex conjugate primes of $\OO_n$.
However, in addition to our example, a rational prime $p$ may or may
not factorise into several primes of $\oo_n$, i.e., we need to look 
at the splitting pattern of rational primes with respect to the double
extension $\KK_{n}/\kk_{n}/\QQ$. Fortunately, the relevant
contribution to the coincidence problem only comes from the so-called
\emph{complex splitting primes}\/ $\mathcal{C}$. They form a subset
of the rational primes $\mathcal{P}$, characterised by splitting
(at least) in the
extension step from $\kk_{n}$ to $\KK_{n}$. Various examples of the
relevant splitting structure are given in \cite{PBR}, see also
\cite{DCG} for their appearance in a related problem.

It was shown in \cite{PBR} that, under {\bf (CN1)}, an integer
\mbox{$k\in\NN$} is a non-trivial coincidence index if and only if it is a
finite product of so-called basic indices \cite{PBR}. They emerge from
the structure of prime ideals in a way that we shall briefly recall
below. The corresponding set of indices is thus a semigroup with unit,
a property that is lost without {\bf (CN1)}.

In our situation with {\bf (CN1)}, the coincidence group is Abelian 
and has the structure
\[
\mathrm{SOC}(\OO_{n}) \;=\; C_{N} \times \ZZ^{(\aleph_{0})},
\]
with $N$ according to Eq.~\eqref{symm}, very much in analogy to
\eqref{3group}, which is the case $n=3$, though with slightly more
complicated expressions for the generators, compare \cite{PBR}.
As above, $\ZZ^{(\aleph_{0})}$ indicates the restriction that
each element is a {\em finite}\/ product of powers of
generators.

\begin{thm} {\rm \bf (CN1)}.
  Let\/ $\mathcal{C}\subset\mathcal{P}$ denote the subset of complex
  splitting primes for the field extension\/ $\KK_{n}/\QQ$. If
  $c^{}_{n}(k)$ denotes the number of simple CSMs of\/ $\OO_{n}$ of index $k$,
  its Dirichlet series generating function is given by
\[
  \varPhi_{n}(s) \; = \; 
  \sum_{k=1}^{\infty} \frac{c^{}_{n}(k)}{k^s} \; = \;
  \prod_{p\in\mathcal{C}} 
  \left(\frac{1+p^{-\ell_p s}}{1-p^{-\ell_p s}}\right)^{m_p/2},
\]
  with characteristic integers\/ $\ell_p$ and\/ $m_p$ as given in
  Tables\/ $\ref{tab:basic}$ and\/ $\ref{tab:rami}$.  This Dirichlet
  series converges absolutely on the half-plane\/
  $\{\mathrm{Re}(s)>1\}$. It can be expressed as a ratio
  of zeta functions,
\[
\varPhi_{n}(s) \;=\; \frac{\zeta^{}_{\KK_{n}}(s)}{\zeta^{}_{\kk_{n}}(2s)}
\cdot
\begin{cases}
(1+p^{-s})^{-1}, & \text{if\/ $n=p^r$,} \\
1, & \text{otherwise,}
\end{cases}
\] 
  where\/ $\zeta^{}_{\FF}(s)$ denotes the Dedekind zeta function of
  the field\/ $\FF$. The first terms of these Dirichlet series
  generating functions are given in Table~$\ref{diritab1}$.
\qed
\label{thm1}
\end{thm}

\noindent {\sc Proof}.
Though a proof was given in \cite{PBR}, we briefly recall the relevant
ingredients for the convenience of the reader. 

A rotation $z=e^{i\varphi}$ produces a CSM of $\OO_n$ if and only if it
can be written as the quotient $u/v$ of two integers $u,v \in \OO_n$,
meaning that $v$ is rotated into $u$, and one has
\[
   \mathrm{SOC}(\OO_n) \; = \;
   \big\{z\in\KK_n \;\bigr|\; \lvert z\rvert=1\big\} .
\]
Moreover, due to {\bf (CN1)}, the integers $u$ and $v$ can always be
chosen coprime (up to units from $\OO_n$), and the representation 
$z=u/v$ is unique in this sense. 

Observe that the \emph{relative norm}\/ of $z\in\KK_n$ over
the maximal real subfield $\kk_n$ is $\mathrm{norm}^{}_{\KK_n/\kk_n}
(z) = \lvert z \rvert^2$, where relative norms of integers (resp.\
units) in $\OO_n$ are integers (resp.\ units) in $\oo_n$. For $z=u/v
\in\mathrm{SOC}(\OO_n)$, with $u,v$ coprime, we clearly must have
\[
    \mathrm{norm}^{}_{\KK_n/\kk_n} (u)  \; = \;
    \mathrm{norm}^{}_{\KK_n/\kk_n} (v)  \; = \;
    \nu \; \in \; \oo_n .
\]
Consequently, every prime factor of $\nu$ in $\oo_n$ must now factor
into two non-associated primes of $\OO_n$, one of which divides $u$
only, while the other divides $v$ only. This restricts us to primes
$\mathfrak{p}\in\oo_n$ that split as $\mathfrak{p} =
\omega_{\mathfrak{p}}\bar{\omega}_{\mathfrak{p}}$ in the extension 
from $\kk_n$ to $\KK_n$, with $\omega_{\mathfrak{p}}/
\bar{\omega}_{\mathfrak{p}}$ {\em not}\/ a unit in $\OO_n$.
They are precisely the ones lying over the rational complex
splitting primes $p$ for the extension $\KK_{n}/\QQ$. All other primes
have cancelled out in the quotient $u/v$, possibly up to a unit. Each 
$z\in\mathrm{SOC}(\OO_n)$ can thus be written as a finite product,
\[
    z \; = \; \varepsilon \prod_{\mathfrak{p}} 
    \left(\frac{\omega^{}_{\mathfrak{p}}}
    {\bar{\omega}^{}_{\mathfrak{p}}}\right)^{t_{\mathfrak{p}}},
\]
where $\mathfrak{p}$ runs over the primes of $\oo_n$ that divide $\nu$
and split as described, and $\varepsilon$ is a unit in $\OO_{n}$.

Given a complex splitting prime $p\in\mathcal{C}$, there are $m_{p}/2$
pairwise non-associated primes $\mathfrak{p}_i \in \oo_n$ that lie
over $p$. Consequently, twice as many primes appear in $\OO_n$, coming
in complex conjugate pairs $\{ \omega^{}_{\mathfrak{p}_i},
\bar{\omega}^{}_{\mathfrak{p}_i} \}$.  Each single such prime has
norm $p^{\ell_p}$ (meaning that the principal ideal defined by it has
index $p^{\ell_p}$ in $\OO_n$).  Moreover, the CSM obtained from a
rotation $z$ as above is the principal ideal generated by its
numerator, which has index
\[
   \prod_{\mathfrak{p}} \big(\mathrm{norm}^{}_{\KK_n/\QQ}
    (\omega^{}_{\mathfrak{p}})\big)^{\lvert t_{\mathfrak{p}}\rvert}
   \; = \; \prod_{\mathfrak{p}} p^{\ell_p \lvert t_{\mathfrak{p}}\rvert}.
\]

Collect now all possibilities that result in the same index. For each
single factor, the corresponding rotation can either be clockwise or
anticlockwise, which give different CSMs. This contributes to
the overall generating function a factor of the form
\[
     \frac{1+x}{1-x} \; = \; 1 + 2x + 2 x^2 + 2x^3 + \ldots
\]
with $x$ replaced by $p^{\ell_p}$, which then appears precisely
$m_{p}/2$ times. Taking the product over all complex splitting
primes, the first claim follows.

The convergence statement is standard, as is the representation
by means of zeta functions. The latter will become clearer from our
further discussion, too.
\qed\smallskip

\begin{table}
\caption[]{Contribution to splitting primes by ramified rational primes 
with corresponding integers $\vphantom{p^{\ell_p}}$
$\ell_p$ and $m_p$ for $\ZZ[\xi_{n}]$, with 
$n$ from the list \eqref{liste}. Here, $r$ is the $p$-free part
$\vphantom{p^{\ell_p}}$ of $n$, 
and $\ell_p\, m_p =\phi(r)$. The corresponding basic index for the 
coincidence spectrum is always $p^{\ell_p}$. \label{tab:rami}}
\def\Strut{\large\strut}
\small
\begin{tabular*}{\linewidth}{@{\extracolsep{\fill}} @{} rrrrrrr @{}}
\hline
$\phi(n)$ &
$n$ &
$p$ &
$r$ &
$\phi(r)$ &
$\ell_p$ &
$m_p$ \Strut\\
\hline
 8 & 20 &  5 & 4 & 2 &  1 & 2 \Strut\\ 
   & 24 &  3 & 8 & 4 &  2 & 2 \Strut\\ \hline
12 & 21 &  7 & 3 & 2 &  1 & 2 \Strut\\ 
   & 28 &  2 & 7 & 6 &  3 & 2 \Strut\\ \hline
16 & 40 &  5 & 8 & 4 &  2 & 2 \Strut\\
   & 48 &  3 &16 & 8 &  4 & 2 \Strut\\ 
   & 60 &  2 &15 & 8 &  4 & 2 \Strut\\
   &    &  3 &20 & 8 &  4 & 2 \Strut\\ 
   &    &  5 &12 & 4 &  2 & 2 \Strut\\ \hline
20 & 33 &  3 &11 &10 &  5 & 2 \Strut\\ \hline
24 & 84 &  2 &21 &12 &  6 & 2 \Strut\\
   &    &  7 &12 & 4 &  2 & 2 \Strut\\
\hline
\end{tabular*}
\end{table}

The integers $\ell_p$ and $m_p$ are characteristic quantities that
depend on $n$ and on the prime $p$. In particular, $m_p$ is the number
of prime ideal divisors of $p$ in $\OO_{n}$, and $\ell_p$ is their
degree (also called the residue class degree) \cite{W}. With the usual
approach of algebraic number theory, its determination (and that of
the corresponding $m_p$) depends on whether $p$ divides $n$ or not.
If $p\nmid n$ , $\ell_p$ is the smallest integer
such that $p^{\ell_p}\equiv 1\; (n)$, and $m_p=\phi(n)/\ell_p$.
Clearly, the
result only depends on the residue class of $p$ modulo $n$.
The numbers $\ell_p$ and $m_p$ are contained in Table~\ref{tab:basic}. 
If, however, $p\mid
n$, one writes $n=p^t r$, where $r$ is the $p$-free part of $n$. Now,
$\ell_p$ is the smallest integer such that $p^{\ell_p}\equiv 1\; (r)$
and $m_p=\phi(r)/\ell_p$, see Table~\ref{tab:rami}. For a step-by-step
derivation of these claims, see Facts 1--3 and Appendix B of
\cite{PBR}.

Referring back to \cite{BG2,W}, and to the theory of Dirichlet characters
summarized there, one can express the zeta functions in terms of
$L$-series as follows,
\begin{align}
  \zeta^{}_{\KK_{n}}(s) & \;=  
      \prod_{\chi\in\widehat{G}_n} L(s,\chi) , 
      \label{L1} \\
  \zeta^{}_{\kk_{n}}(s) & \; = 
      \prod_{\substack{\chi\in\widehat{G}_n \\ 
      \chi \text{ even}}} L(s,\chi),
      \label{L2}
\end{align}
where a character $\chi$ is called {\em even}\/ if $\chi(-1)=1$ and the
$L$-series read
\[
   L(s,\chi) \; = \; \sum_{k=1}^{\infty} 
   \frac{\chi(k)}{k^s} .
\]
As in \cite{BG2,W}, $\widehat{G}_n$ denotes the set of primitive
Dirichlet characters for the field extension $\KK_n/\QQ$, and
thus contains $\phi(n)$ elements. The principal character in this
formulation, $\chi^{}_{0} \equiv 1$, always leads to
$L(s,\chi^{}_{0})=\zeta(s)$, i.e., to Riemann's zeta function itself.
This is a meromorphic function on $\CC$, with a single simple pole
of residue $1$ at $s=1$, and no zeros in the closed half-plane
$\{s\in\CC\mid \mathrm{Re}(s)\ge 1\}$.
All remaining $L$-series are entire functions
on $\CC$, without any zeros in this half-plane either, compare
\cite[Thm.~12.5]{A} or \cite[Ch.~4]{W}.

Note that
\[
    \zeta^{}_{\KK_n}(s) \; = \; 
    \sum_{k=1}^{\infty} \frac{a_n (k)}{k^s},
\]
under the assumption ({\bf CN1}), counts the {\em principal}\/ ideals of
$\OO_n$, so that $a_n (k)$ is the number of similarity sublattices
or submodules of $\OO_n$ of index $k$. It was shown in \cite[Thm.~2]{BG2} that 
\[
    A_{n}(x) \; := \;  \sum_{k\le x} a_{n}(k) \; \sim \; \alpha^{}_{n} x ,
\]
with the growth constant
\[
    \alpha^{}_{n} \; = \; \mathrm{res}^{}_{s=1}\big(\zeta^{}_{\KK_{n}}(s)\big)
    \; = \prod_{1\not\equiv\chi\in\widehat{G}_{n}} L(1,\chi).
\]
Examples are given in Table~\ref{tab:asym}.

Let us next explain how to calculate the coefficients $c^{}_{n}(k)$
defined in Theorem~\ref{thm1} explicitly.  Under our assumption ({\bf
CN1}), the arithmetic function $c^{}_{n}(k)$ is multiplicative, whence
it is sufficient to know the values $c^{}_{n}(p^r)$ for all primes $p$
and powers $r>0$. On the level of the generating function, this
corresponds to expanding the factors $E(p^{-s})$ of the Euler products
$\prod_{p} E (p^{-s})$, where $p$ runs over the primes of $\ZZ$.

Here and below, each Euler factor is either of the form
\begin{align}
   E (p^{-s}) & \; = \;  \frac{1}{(1-p^{-\ell s})^m}
   \label{factor1} \\
 & \; = \; 
 \sum_{j=0}^{\infty} \binom{j+m-1}{m-1}
 \frac{1}{(p^{s})^{\ell j}}  \notag
\end{align}
or, as for $\varPhi(s)$ above, reads
\begin{align}
   E (p^{-s}) & \; = \;  \left(\frac{1+p^{-\ell s}}{1-p^{-\ell s}}\right)^{m} 
   \label{factor2} \\
 & \; = \; 
 \sum_{j=0}^{\infty} \left(\sum_{i=0}^{j}
 \binom{m}{i} \binom{m-1+j-i}{m-1}\right)
 \frac{1}{(p^{s})^{\ell j}}.  \notag
\end{align}
These two formulas follow easily from expanding 
\[
   \left(\frac{1}{1-x}\right)^m \quad \mbox{ and } \quad
   \left(\frac{1+x}{1-x}\right)^m
\]
into power series, followed by inserting $p^{-\ell s}$ for $x$.
{}From \eqref{factor1} and \eqref{factor2}, one then quickly extracts the 
values of $c^{}_{n}(p^r)$ for $r\ge 0$.

\begin{table}[t]
\caption[]{Numerical values of residues (at $s=1$) of Dedekind zeta functions
and related Dirichlet series generating functions.
\label{tab:asym}}\smallskip
\def\Strut{\large\strut}
\small
\renewcommand{\arraystretch}{1}
\begin{tabular*}{\linewidth}{@{\extracolsep{\fill}} @{} 
cr@{\qquad}c@{\qquad}c@{\qquad}c @{}}
\hline
$\phi(n)$ & $n$ & $\alpha^{}_{n}$ & $\beta^{}_{n}$ & $\gamma^{}_{n}$ \Strut\\
\hline
2 &  3 & $0.604\,600$ & $0.285\,041$ & $0.275\,664$ \rule[0ex]{0ex}{3ex}\Strut\\
  &  4 & $0.785\,398$ & $0.336\,193$ & $0.318\,310$ \rule[-1.5ex]{0ex}{1.5ex}\Strut\\
\hline
4 &  5 & $0.339\,837$ & $0.249\,136$ & $0.243\,785$ \rule[0ex]{0ex}{3ex}\Strut\\
  &  8 & $0.543\,676$ & $0.258\,663$ & $0.252\,584$ \Strut\\
  & 12 & $0.361\,051$ & $0.235\,129$ & $0.231\,117$ \rule[-1.5ex]{0ex}{1.5ex}\Strut\\
\hline
6 &  7 & $0.287\,251$ & $0.240\,733$ & $0.235\,393$ \rule[0ex]{0ex}{3ex}\Strut\\
  &  9 & $0.333\,685$ & $0.216\,086$ & $0.213\,490$ \rule[-1.5ex]{0ex}{1.5ex}\Strut\\
\hline
8 & 15 & $0.215\,279$ & $0.202\,672$ & $0.200\,651$ \rule[0ex]{0ex}{3ex}\Strut\\
  & 16 & $0.464\,557$ & $0.230\,565$ & $0.226\,884$ \Strut\\
  & 20 & $0.288\,769$ & $0.267\,051$ & $0.255\,192$ \Strut\\
  & 24 & $0.299\,995$ & $0.222\,573$ & $0.218\,534$ \rule[-1.5ex]{0ex}{1.5ex}\Strut\\
\hline
10& 11 & $0.239\,901$ & $0.217\,016$ & $0.214\,504$ \rule[-1.5ex]{0ex}{4.5ex}\Strut\\
\hline
12& 13 & $0.213\,514$ & $0.195\,868$ & $0.194\,563$ \rule[0ex]{0ex}{3ex}\Strut\\
  & 21 & $0.227\,271$ & $0.225\,794$ & $0.220\,216$ \Strut\\
  & 28 & $0.251\,795$ & $0.245\,523$ & $0.239\,718$ \Strut\\
  & 36 & $0.220\,933$ & $0.192\,881$ & $0.191\,637$ \rule[-1.5ex]{0ex}{1.5ex}\Strut\\
\hline
16 & 17 & $0.153\,708$ & $0.143\,107$ & $0.142\,875$ \rule[0ex]{0ex}{3ex}\Strut\\
   & 32 & $0.278\,432$ & $0.137\,913$ & $0.137\,709$ \Strut\\
   & 40 & $0.210\,946$ & $0.197\,390$ & $0.195\,931$ \Strut\\
   & 48 & $0.222\,767$ & $0.166\,373$ & $0.165\,820$ \Strut\\
   & 60 & $0.195\,032$ & $0.194\,455$ & $0.192\,773$ \rule[-1.5ex]{0ex}{1.5ex}\Strut\\
\hline
18 & 19 &  $0.121\,018$ & $0.113\,726$ & $0.113\,651$ \rule[0ex]{0ex}{3ex}\Strut\\
   & 27 &  $0.212\,854$ & $0.141\,292$ & $0.141\,081$ \rule[-1.5ex]{0ex}{1.5ex}\Strut\\ 
\hline
20 & 25 &  $0.181\,458$ & $0.144\,724$ & $0.144\,466$ \rule[0ex]{0ex}{3ex}\Strut\\
   & 33 &  $0.159\,226$ & $0.157\,420$ & $0.156\,974$ \Strut\\ 
   & 44 &  $0.126\,912$ & $0.124\,934$ & $0.124\,741$ \rule[-1.5ex]{0ex}{1.5ex}\Strut\\
\hline
24 & 35 &  $0.166\,239$ & $0.166\,011$ & $0.165\,506$ \rule[0ex]{0ex}{3ex}\Strut\\
   & 45 &  $0.118\,121$ & $0.116\,387$ & $0.116\,302$ \Strut\\
   & 84 &  $0.116\,090$ & $0.115\,471$ & $0.115\,336$ \rule[-1.5ex]{0ex}{1.5ex}\Strut\\
\hline
\end{tabular*}
\renewcommand{\arraystretch}{1}
\end{table}

The Dirichlet series generating function of an arithmetic function permits 
the derivation of asymptotic properties (here, of $c^{}_{n}(k)$) or, 
more precisely, of the corresponding summatory function, 
$\sum_{k\le x} c^{}_{n}(k)$, as
$x\to\infty$. Here, Delange's theorem (see \cite[Ch.~II.7, Thm.~15]{T}
for the general formulation, and \cite[Prop.~4]{BG2} for the reduction
to the situation at hand) leads to the following result.

\begin{coro}
{\rm \bf (CN1)}. 
The number of simple CSMs of index $\le x$ is given by
\[
   \sum_{k\le x}c^{}_{n}(k) \; \sim \;
   x\cdot \big(\mbox{\rm res}_{s=1} \,\varPhi_{n}(s)\big)
   \; = \; x\cdot \gamma_n ,
\]
with the residue
\[
   \gamma_n \; = \;  \frac{\alpha^{}_{n}}
   {\zeta^{}_{\kk_{n}}(2)} \cdot
   \begin{cases}
          p/(p+1), & \text{if\/ $n=p^r$,} \\
          1, & \text{otherwise.}
   \end{cases}
\]
Here, $\alpha^{}_{n} = \mbox{\rm res}^{}_{s=1}\,\zeta^{}_{\KK_{n}}(s)=
\prod_{1\not\equiv\chi\in\widehat{G}_{n}}L(1,\chi)$, with\/
$\widehat{G}_{n}$ denoting the set of primitive Dirichlet
characters of\/ $\KK_n/\QQ$.
\label{coro1}
\end{coro}
\noindent {\sc Proof}.
The coefficients $c^{}_n (k)$ are non-negative numbers, and their
generating function $\varPhi_n (s)$, by Theorem~\ref{thm1}, is
a meromorphic function on the half-plane $\{ \mathrm{Re}(s)> 1/2\}$,
with a single simple pole at $s=1$. Delange's theorem then results in
the asymptotic linear growth as claimed, with the residue of
$\varPhi_n (s)$ at $s=1$ as growth constant. Its calculation is obvious
from Theorem~\ref{thm1} together with formulas \eqref{L1} and \eqref{L2}.
\qed\smallskip

One way to calculate the residues uses the factorization of our zeta 
functions into $L$-series, followed by rewriting the latter in terms 
of the Hurwitz zeta function \cite[Ch.~12.1]{A}
\begin{equation}
  \zeta(s,a) \; = \; \sum_{k=0}^{\infty}
  \frac{1}{(k+a)^s} , \label{Hurwitz1} 
\end{equation}
which is well defined for $0<a\le 1$ and absolutely convergent on the
half-plane $\{\mathrm{Re}(s)>1\}$. On $\CC$, analytic continuation
defines $\zeta(s,a)$ as a meromorphic function with a single simple
pole of residue $1$ at $s=1$, see \cite[Thm.~12.4]{A}. In particular,
$\zeta(s,1)=\zeta(s)$. Furthermore, if $\Gamma (s)$ denotes the
ordinary gamma function, one has
\begin{equation}
  \lim_{s\to 1} \Big(\zeta(s,a)-\frac{1}{s-1}\Big) \; = \;
  - \frac{\Gamma'(a)}{\Gamma(a)},
\label{gamma}
\end{equation}
see \cite[p.~271]{WW}. For further connections to the generalised
Euler constants and the Laurent series of the Hurwitz zeta function
at $s=1$, see \cite{Berndt}. 

Next, if a Dirichlet character has
period $q$\/ (e.g., if $q=f_\chi$ is the conductor of $\chi$, which is
its minimal period), one has the representation \cite{A}
\begin{equation}
  L(s,\chi) \; = \; q^{-s} \sum_{r=1}^{q} \chi(r)\, \zeta(s,\tfrac{r}{q}).
\label{Hurwitz2}
\end{equation}
Moreover, if $\chi$ is not the principal character, one can rewrite this
equation as
\begin{equation}
  L(s,\chi) \; = \; q^{-s} \sum_{r=1}^{q} \chi(r)\, 
  \Big(\zeta(s,\tfrac{r}{q})-\frac{1}{s-1}\Big),
\label{Hurwitz3}
\end{equation}
because $\sum_{r=1}^{q}\chi(r)=0$ for non-principal characters by the
orthogonality relation. This shows regularity at $s=1$, and with 
Eq.~\eqref{gamma} one obtains
\begin{equation}
  L(1,\chi) \; = \; - \frac{1}{q} \sum_{r=1}^{q} \chi(r)\, 
\frac{\Gamma'(\tfrac{r}{q})}{\Gamma(\tfrac{r}{q})}.
\label{Hurwitz4}
\end{equation}

In \cite{BG2}, we additionally used the class number formula to
transform the calculation of $\mathrm{res}_{s=1}
\big(\zeta^{}_{\KK_n}(s)\big)$ into an algebraic problem. However,
given the fact that $\zeta(s,a)$ and $\Gamma(s)$ are accessible by means of
algebraic program packages to arbitrary numerical precision, it seems
easier to directly use formulas
\eqref{L1} and \eqref{L2} together with \eqref{Hurwitz2} for these
calculations. Some results\footnote{Please note that the exact value 
of $\gamma^{}_{7}$ shown in \cite{PBR} should be replaced by 
$\gamma^{}_{7}=\frac{21\sqrt{7}R}{16\pi^3}$, with
the regulator $R$ as given there.} are given in Table~\ref{tab:asym}.

\section{Multiple coincidences}

In generalisation of what we explained in 
Section~\ref{example:triangular}, we are now interested in multiple
intersections of the form
\[
\OO_{n}\cap e^{i\varphi^{}_{1}}\OO_{n}\cap\ldots\cap 
e^{i\varphi^{}_{m}}\OO_{n},
\]
provided that this is a submodule of full rank, and the corresponding index
\[
\Sigma(\varphi^{}_{1},\ldots,\varphi^{}_{m}) \; := \;
[\OO_{n} : (\OO_{n}\cap e^{i\varphi^{}_{1}}\OO_{n}\cap\ldots\cap 
e^{i\varphi^{}_{m}}\OO_{n})].
\]
We call such intersections \emph{multiple coincidence site modules},
or multiple CSMs for short.

\begin{thm} {\rm \bf (CN1)}.
  Let\/  $b_{n}(k)$ denotes the number of multiple CSMs of\/ $\OO_{n}$ 
  of index $k$, including the simple ones. Then,
  its Dirichlet series generating function is given by
\[
   \varPsi_{n}(s) \; = \; 
   \sum_{k=1}^{\infty} \frac{b_{n}(k)}{k^s} \; = \;
   \prod_{p\in\mathcal{C}} 
   \left(\frac{1}{1-p^{-\ell_p s}}\right)^{m_p},
\]
  with\/ $\mathcal C$ once again the complex splitting primes
  of the field extension\/ $\KK_n/\QQ$
  and\/ $\ell_p$, $m_p$ the characteristic integers introduced
  above.
  This Dirichlet series converges absolutely on the half-plane\/
  $\{\mathrm{Re}(s)>1\}$. The first terms of the difference\/
  $\varPsi_n (s)-\varPhi_n (s)$ are given in Table~$\ref{diritab2}$.
\label{thm2}
\end{thm}

\noindent {\sc Proof}.
Observe as before that 
\[
\OO_{n}\cap e^{i\varphi^{}_{1}}\OO_{n}\cap\ldots\cap e^{i\varphi^{}_{m}}\OO_{n}
\; = \; \bigcap_{\ell=1}^{m}\, (\OO_{n}\cap e^{i\varphi^{}_{\ell}}\OO_{n}).
\]
This shows that the multiple intersection is a submodule of $\OO_{n}$
of full rank $\phi(n)$ if and only if each $\OO_{n}\cap
e^{i\varphi^{}_{\ell}}\OO_{n}$ is a simple CSM of $\OO_{n}$, each of
which is a principal ideal of $\OO_{n}$ due to {\bf (CN1)}. Consequently,
the multiple intersection is the lcm of the simple CSMs involved, as in
Eq.~\eqref{lcm-eq}.

For each prime ideal $\mathfrak{p}$ of $\oo_{n}$ that lies over the
complex splitting prime $p$, the combinatorial argument used in the
proof of Proposition~\ref{prop1} applies \emph{independently}. This
modifies each factor of the Euler product decomposition of
Theorem~\ref{thm1} in the same way as in our previous example (the
triangular lattice), i.e., each term of the form $(1+x)$ in the
numerator is replaced by a term $(1-x)^{-1}$, with $x=p^{-\ell_{p}s}$.
This gives the new Euler product, which clearly converges as claimed.
\qed\smallskip

\begin{table*}
\caption[]{First terms of the Dirichlet series from Theorem~\ref{thm1}
for simple coincidences, for integers from the list \eqref{liste}; 
see also \cite{GP}.\label{diritab1}}
\def\Strut{\large\strut}
\renewcommand{\arraystretch}{2.2}
\small
\begin{tabular*}{\textwidth}{@{\extracolsep{\fill}} @{} rl @{}}
\hline
 $n$ & $\varPhi_{n}(s)$\Strut\\
\hline
 3 & $1 + \frac{2}{7^s}\+\frac{2}{13^s}\+\frac{2}{19^s}\+\frac{2}{31^s}\+ 
\frac{2}{37^s}\+\frac{2}{43^s}\+\frac{2}{49^s}\+\frac{2}{61^s}\+\frac{2}{67^s}
\+\frac{2}{73^s}\+\frac{2}{79^s}\+\frac{4}{91^s}\+\frac{2}{97^s}\+ 
\frac{2}{103^s}\+\frac{2}{109^s}\+\frac{2}{127^s}\+\frac{4}{133^s}
+ \ldots$\Strut\\
 4 & $1 + \frac{2}{5^s}\+ \frac{2}{13^s}\+ \frac{2}{17^s}\+ \frac{2}{25^s}\+
 \frac{2}{29^s}
\+ \frac{2}{37^s}\+ \frac{2}{41^s}\+ \frac{2}{53^s}\+ \frac{2}{61^s}\+ 
\frac{4}{65^s}\+ \frac{2}{73^s}\+ \frac{4}{85^s}\+ \frac{2}{89^s}\+ 
\frac{2}{97^s}\+ \frac{2}{101^s}\+ \frac{2}{109^s}\+ \frac{2}{113^s}
+ \ldots$\Strut\\
 5 & $1 + \frac{4}{11^s}\+\frac{4}{31^s}\+\frac{4}{41^s}\+\frac{4}{61^s}
\+\frac{4}{71^s}\+\frac{4}{101^s}\+\frac{8}{121^s}\+\frac{4}{131^s}\+\frac{4}{151^s}
\+\frac{4}{181^s}\+\frac{4}{191^s}\+\frac{4}{211^s}\+\frac{4}{241^s}\+\frac{4}{251^s}
\+\frac{4}{271^s}\+\frac{4}{281^s}
+ \ldots$\Strut\\
 7 & $1 + \frac{2}{8^s}\+\frac{6}{29^s}\+\frac{6}{43^s}\+\frac{2}{64^s}\+\frac{6}{71^s}
\+\frac{6}{113^s}\+\frac{6}{127^s}\+\frac{6}{197^s}\+\frac{6}{211^s}\+\frac{12}{232^s}
\+\frac{6}{239^s}\+\frac{6}{281^s}\+\frac{6}{337^s}\+\frac{12}{344^s}\+\frac{6}{379^s}
\+\frac{6}{421^s}
+ \ldots $\Strut\\
 8 & $1 + \frac{2}{9^s}\+\frac{4}{17^s}\+\frac{2}{25^s}\+\frac{4}{41^s}
\+\frac{4}{73^s}\+\frac{2}{81^s}\+\frac{4}{89^s}\+\frac{4}{97^s}\+\frac{4}{113^s}
\+\frac{2}{121^s}\+\frac{4}{137^s}\+\frac{8}{153^s}\+\frac{2}{169^s}\+\frac{4}{193^s}
\+\frac{4}{225^s}\+\frac{4}{233^s}
+ \ldots$\Strut\\
 9 & $1 + \frac{6}{19^s}\+\frac{6}{37^s}\+\frac{6}{73^s}\+\frac{6}{109^s}
\+\frac{6}{127^s}\+\frac{6}{163^s}\+\frac{6}{181^s}\+\frac{6}{199^s}\+\frac{6}{271^s}
\+\frac{6}{307^s}\+\frac{2}{343^s}\+\frac{18}{361^s}\+\frac{6}{379^s}\+\frac{6}{397^s}
\+\frac{6}{433^s}\+\frac{6}{487^s}
+ \ldots $\Strut\\
11 & $1 + \frac{10}{23^s}\+\frac{10}{67^s}\+\frac{10}{89^s}\+\frac{10}{199^s}
\+\frac{2}{243^s}\+\frac{10}{331^s}\+\frac{10}{353^s}\+\frac{10}{397^s}
\+\frac{10}{419^s}\+\frac{10}{463^s}\+\frac{50}{529^s}\+\frac{10}{617^s}
\+\frac{10}{661^s}\+\frac{10}{683^s}\+\frac{10}{727^s}\+\frac{10}{859^s}
+ \ldots$\Strut\\
12 & $1 + \frac{4}{13^s}\+\frac{2}{25^s}\+\frac{4}{37^s}\+\frac{2}{49^s}
\+\frac{4}{61^s}\+\frac{4}{73^s}\+\frac{4}{97^s}\+\frac{4}{109^s}
\+\frac{4}{157^s}\+\frac{8}{169^s}\+\frac{4}{181^s}\+\frac{4}{193^s}
\+\frac{4}{229^s}\+\frac{4}{241^s}\+\frac{4}{277^s}\+\frac{2}{289^s}
 +\ldots$\Strut\\
13 & $1 + \frac{4}{27^s}\+\frac{12}{53^s}\+\frac{12}{79^s}\+\frac{12}{131^s}
\+\frac{12}{157^s}\+\frac{12}{313^s}\+\frac{12}{443^s}\+\frac{12}{521^s}
\+\frac{12}{547^s}\+\frac{12}{599^s}\+\frac{12}{677^s}\+\frac{8}{729^s}
\+\frac{12}{859^s}\+\frac{12}{911^s}\+\frac{12}{937^s}
+ \ldots$\Strut\\
15 & $1 + \frac{2}{16^s}\+\frac{8}{31^s}\+\frac{8}{61^s}\+\frac{4}{121^s}
\+\frac{8}{151^s}\+\frac{8}{181^s}\+\frac{8}{211^s}\+\frac{8}{241^s}
\+\frac{2}{256^s}\+\frac{8}{271^s}\+\frac{8}{331^s}\+\frac{4}{361^s}
\+\frac{8}{421^s}\+\frac{16}{496^s}\+\frac{8}{541^s}\+\frac{8}{571^s}
+ \ldots$\Strut\\
16 & $1 + \frac{8}{17^s}\+\frac{4}{49^s}\+\frac{2}{81^s}\+\frac{8}{97^s}
\+\frac{8}{113^s}\+\frac{8}{193^s}\+\frac{8}{241^s}\+\frac{8}{257^s}
\+\frac{32}{289^s}\+\frac{8}{337^s}\+\frac{8}{353^s}\+\frac{8}{401^s}
\+\frac{8}{433^s}\+\frac{8}{449^s}\+\frac{4}{529^s}\+\frac{8}{577^s}
+ \ldots$\Strut\\
17 & $1 + \frac{16}{103^s}\+\frac{16}{137^s}\+\frac{16}{239^s}
\+\frac{16}{307^s}\+\frac{16}{409^s}\+\frac{16}{443^s}\+\frac{16}{613^s}
\+\frac{16}{647^s}\+\frac{16}{919^s}\+\frac{16}{953^s}\+\frac{16}{1021^s}
\+\frac{16}{1123^s}\+\frac{16}{1259^s}\+\frac{16}{1327^s}\+\frac{16}{1361^s}
+\ldots$\Strut\\
19 & $1 + \frac{18}{191^s}\+\frac{18}{229^s}\+\frac{6}{343^s}
\+\frac{18}{419^s}\+\frac{18}{457^s}\+\frac{18}{571^s}\+\frac{18}{647^s}
\+\frac{18}{761^s}\+\frac{18}{1103^s}\+\frac{18}{1217^s}\+\frac{6}{1331^s}
\+\frac{18}{1483^s}\+\frac{18}{1559^s}\+\frac{18}{1597^s}
+ \ldots$\Strut\\
20 & $1 + \frac{2}{5^s}\+\frac{2}{25^s}\+\frac{8}{41^s}\+\frac{8}{61^s}
\+\frac{2}{81^s}\+\frac{8}{101^s}\+\frac{4}{121^s}\+\frac{2}{125^s}
\+\frac{8}{181^s}\+\frac{16}{205^s}\+\frac{8}{241^s}\+\frac{8}{281^s}
\+\frac{16}{305^s}\+\frac{8}{401^s}\+\frac{4}{405^s}\+\frac{8}{421^s}
+ \ldots$\Strut\\
21 & $1 + \frac{2}{7^s}\+\frac{12}{43^s}\+\frac{2}{49^s}\+\frac{2}{64^s}
\+\frac{12}{127^s}\+\frac{6}{169^s}\+\frac{12}{211^s}\+\frac{24}{301^s}
\+\frac{12}{337^s}\+\frac{2}{343^s}\+\frac{12}{379^s}\+\frac{12}{421^s}
\+\frac{4}{448^s}\+\frac{12}{463^s}\+\frac{12}{547^s}\+\frac{12}{631^s}
+ \ldots$\Strut\\
24 & $1 + \frac{2}{9^s}\+\frac{4}{25^s}\+\frac{4}{49^s}\+\frac{8}{73^s}
\+\frac{2}{81^s}\+\frac{8}{97^s}\+\frac{4}{121^s}\+\frac{4}{169^s}
\+\frac{8}{193^s}\+\frac{8}{225^s}\+\frac{8}{241^s}\+\frac{4}{289^s}
\+\frac{8}{313^s}\+\frac{8}{337^s}\+\frac{4}{361^s}\+\frac{8}{409^s}
+ \ldots$\Strut\\
25 & $1 + \frac{20}{101^s}\+\frac{20}{151^s}\+\frac{20}{251^s}
\+\frac{20}{401^s}\+\frac{20}{601^s}\+\frac{20}{701^s}\+\frac{20}{751^s}
\+\frac{20}{1051^s}\+\frac{20}{1151^s}\+\frac{20}{1201^s}
\+\frac{20}{1301^s}\+\frac{20}{1451^s}\+\frac{20}{1601^s}
\+\frac{20}{1801^s}
+ \ldots$\Strut\\
27 & $1 + \frac{18}{109^s}\+\frac{18}{163^s}\+\frac{18}{271^s}
\+\frac{18}{379^s}\+\frac{18}{433^s}\+\frac{18}{487^s}
\+\frac{18}{541^s}\+\frac{18}{757^s}\+\frac{18}{811^s}
\+\frac{18}{919^s}\+\frac{18}{1297^s}\+\frac{18}{1459^s}
\+\frac{18}{1567^s}\+\frac{18}{1621^s}\+\frac{18}{1783^s}
+ \ldots$\Strut\\
28 & $1 + \frac{2}{8^s}\+\frac{12}{29^s}\+\frac{2}{64^s}
\+\frac{12}{113^s}\+\frac{6}{169^s}\+\frac{12}{197^s}
\+\frac{24}{232^s}\+\frac{12}{281^s}\+\frac{12}{337^s}
\+\frac{12}{421^s}\+\frac{12}{449^s}\+\frac{2}{512^s}
\+\frac{12}{617^s}\+\frac{12}{673^s}\+\frac{12}{701^s}
\+\frac{12}{757^s}
+ \ldots$\Strut\\
32 & $1 + \frac{16}{97^s}\+\frac{16}{193^s}\+\frac{16}{257^s}
\+\frac{8}{289^s}\+\frac{16}{353^s}\+\frac{16}{449^s}
\+\frac{16}{577^s}\+\frac{16}{641^s}\+\frac{16}{673^s}
\+\frac{16}{769^s}\+\frac{16}{929^s}\+\frac{16}{1153^s}
\+\frac{16}{1217^s}\+\frac{16}{1249^s}\+\frac{16}{1409^s}
+ \ldots$\Strut\\
33 & $1 + \frac{20}{67^s}\+\frac{20}{199^s}\+\frac{2}{243^s}
\+\frac{20}{331^s}\+\frac{20}{397^s}\+\frac{20}{463^s}
\+\frac{10}{529^s}\+\frac{20}{661^s}\+\frac{20}{727^s}
\+\frac{20}{859^s}\+\frac{20}{991^s}\+\frac{20}{1123^s}
\+\frac{20}{1321^s}\+\frac{20}{1453^s}\+\frac{20}{1783^s}
+ \ldots$\Strut\\
35 & $1 + \frac{24}{71^s}\+\frac{24}{211^s}\+\frac{24}{281^s}
\+\frac{24}{421^s}\+\frac{24}{491^s}\+\frac{24}{631^s}
\+\frac{24}{701^s}\+\frac{12}{841^s}\+\frac{24}{911^s}
\+\frac{24}{1051^s}\+\frac{8}{1331^s}\+\frac{24}{1471^s}
\+\frac{12}{1681^s}\+\frac{24}{2311^s}\+\frac{24}{2381^s}
+ \ldots$\Strut\\
36 & $1 + \frac{12}{37^s}\+\frac{12}{73^s}\+\frac{12}{109^s}
\+\frac{12}{181^s}\+\frac{6}{289^s}\+\frac{6}{361^s}
\+\frac{12}{397^s}\+\frac{12}{433^s}\+\frac{12}{541^s}
\+\frac{12}{577^s}\+\frac{12}{613^s}\+\frac{12}{757^s}
\+\frac{12}{829^s}\+\frac{12}{937^s}\+\frac{12}{1009^s}
+ \ldots$\Strut\\
40 & $1 + \frac{2}{25^s}\+\frac{16}{41^s}\+\frac{4}{81^s}
\+\frac{8}{121^s}\+\frac{16}{241^s}\+\frac{16}{281^s}
\+\frac{8}{361^s}\+\frac{16}{401^s}\+\frac{16}{521^s}
\+\frac{16}{601^s}\+\frac{2}{625^s}\+\frac{16}{641^s}
\+\frac{16}{761^s}\+\frac{8}{841^s}\+\frac{16}{881^s}
+ \ldots$\Strut\\
44 & $1 + \frac{20}{89^s}\+\frac{20}{353^s}\+\frac{20}{397^s}
\+\frac{10}{529^s}\+\frac{20}{617^s}\+\frac{20}{661^s}
\+\frac{20}{881^s}\+\frac{20}{1013^s}\+\frac{20}{1277^s}
\+\frac{20}{1321^s}\+\frac{20}{1409^s}\+\frac{20}{1453^s}
\+\frac{20}{2069^s}\+\frac{20}{2113^s}
+ \ldots$\Strut\\
45 & $1 + \frac{24}{181^s}\+\frac{24}{271^s}\+\frac{12}{361^s}
\+\frac{24}{541^s}\+\frac{24}{631^s}\+\frac{24}{811^s}
\+\frac{24}{991^s}\+\frac{24}{1171^s}\+\frac{24}{1531^s}
\+\frac{24}{1621^s}\+\frac{24}{1801^s}\+\frac{24}{2161^s}
\+\frac{24}{2251^s}\+\frac{24}{2341^s}
+  \ldots$\Strut\\
48 & $1 + \frac{8}{49^s}\+\frac{2}{81^s}\+\frac{16}{97^s}
\+\frac{16}{193^s}\+\frac{16}{241^s}\+\frac{8}{289^s}
\+\frac{16}{337^s}\+\frac{16}{433^s}\+\frac{8}{529^s}
\+\frac{16}{577^s}\+\frac{4}{625^s}\+\frac{16}{673^s}
\+\frac{16}{769^s}\+\frac{8}{961^s}\+\frac{16}{1009^s}
+ \ldots$\Strut\\
60 & $1 + \frac{2}{16^s}\+\frac{2}{25^s}\+\frac{16}{61^s}
\+\frac{2}{81^s}\+\frac{8}{121^s}\+\frac{16}{181^s}
\+\frac{16}{241^s}\+\frac{2}{256^s}\+\frac{8}{361^s}
\+\frac{4}{400^s}\+\frac{16}{421^s}\+\frac{16}{541^s}
\+\frac{16}{601^s}\+\frac{2}{625^s}\+\frac{16}{661^s}
\+\frac{8}{841^s}
+ \ldots$\Strut\\
84 & $1 + \frac{2}{49^s}\+\frac{2}{64^s}\+\frac{12}{169^s}
\+\frac{24}{337^s}\+\frac{24}{421^s}\+\frac{24}{673^s}
\+\frac{24}{757^s}\+\frac{12}{841^s}\+\frac{24}{1009^s}
\+\frac{24}{1093^s}\+\frac{24}{1429^s}\+\frac{24}{1597^s}
\+\frac{12}{1681^s}\+\frac{12}{1849^s}\+\frac{24}{1933^s}
+ \ldots$\Strut\\[2mm] 
\hline
\end{tabular*}
\renewcommand{\arraystretch}{1}
\end{table*}

\begin{table*}
\caption[]{Additional terms in the Dirichlet series of Theorem~\ref{thm2}
for truly multiple coincidences; see also \cite{GP}.\label{diritab2}}
\def\Strut{\large\strut}
\renewcommand{\arraystretch}{2.2}
\small
\begin{tabular*}{\textwidth}{@{\extracolsep{\fill}} @{} rl @{}}
\hline
 $n$ & $\varPsi_{n}(s)-\varPhi_{n}(s)$\Strut\\
\hline
 3 & $\frac{1}{49^s}\+\frac{1}{169^s}\+\frac{2}{343^s}\+\frac{1}{361^s}
\+\frac{2}{637^s}\+\frac{2}{931^s}\+\frac{1}{961^s}\+\frac{2}{1183^s}
\+\frac{1}{1369^s}\+\frac{2}{1519^s}\+\frac{2}{1813^s}\+\frac{1}{1849^s}
\+\frac{2}{2107^s}\+\frac{2}{2197^s}\+\frac{3}{2401^s}
+ \ldots$\Strut\\
 4 & $\frac{1}{25^s}\+\frac{2}{125^s}\+\frac{1}{169^s}\+\frac{1}{289^s}
\+\frac{2}{325^s}\+\frac{2}{425^s}\+\frac{3}{625^s}\+\frac{2}{725^s}
\+\frac{1}{841^s}\+\frac{2}{845^s}\+\frac{2}{925^s}\+\frac{2}{1025^s}
\+\frac{2}{1325^s}\+\frac{1}{1369^s}\+\frac{2}{1445^s}
+ \ldots$\Strut\\
 5 & $\frac{2}{121^s}\+\frac{2}{961^s}\+\frac{8}{1331^s}\+\frac{2}{1681^s}
\+\frac{2}{3721^s}\+\frac{8}{3751^s}\+\frac{8}{4961^s}\+\frac{2}{5041^s}
\+\frac{8}{7381^s}\+\frac{8}{8591^s}\+\frac{2}{10201^s}\+\frac{8}{10571^s}
\+\frac{8}{12221^s}\+\frac{19}{14641^s}
+ \ldots$\Strut\\
 7 & $\frac{1}{64^s}\+\frac{2}{512^s}\+\frac{3}{841^s}\+\frac{3}{1849^s}
\+\frac{6}{1856^s}\+\frac{6}{2752^s}\+\frac{3}{4096^s}\+\frac{6}{4544^s}
\+\frac{3}{5041^s}\+\frac{6}{6728^s}\+\frac{6}{7232^s}\+\frac{6}{8128^s}
\+\frac{6}{12608^s}\+\frac{3}{12769^s}
+ \ldots $\Strut\\
 8 & $\frac{1}{81^s}\+\frac{2}{289^s}\+\frac{1}{625^s}\+\frac{2}{729^s}
\+\frac{4}{1377^s}\+\frac{2}{1681^s}\+\frac{2}{2025^s}\+\frac{4}{2601^s}
\+\frac{4}{3321^s}\+\frac{8}{4913^s}\+\frac{2}{5329^s}\+\frac{2}{5625^s}
\+\frac{4}{5913^s}\+\frac{3}{6561^s}
+ \ldots$\Strut\\
 9 & $\frac{3}{361^s}\+\frac{3}{1369^s}\+\frac{3}{5329^s}\+\frac{18}{6859^s}
\+\frac{3}{11881^s}\+\frac{18}{13357^s}\+\frac{3}{16129^s}\+\frac{18}{26011^s}
\+\frac{18}{26353^s}\+\frac{3}{26569^s}\+\frac{3}{32761^s}\+\frac{18}{39349^s}
\+\frac{3}{39601^s}
+ \ldots $\Strut\\
11 & $\frac{5}{529^s}\+\frac{5}{4489^s}\+\frac{5}{7921^s}\+\frac{50}{12167^s}
\+\frac{50}{35443^s}\+\frac{5}{39601^s}\+\frac{50}{47081^s}\+\frac{1}{59049^s}
\+\frac{50}{103247^s}\+\frac{50}{105271^s}\+\frac{5}{109561^s}
\+\frac{5}{124609^s}
+ \ldots$\Strut\\
12 & $\frac{2}{169^s}\+\frac{1}{625^s}\+\frac{2}{1369^s}\+\frac{8}{2197^s}
\+\frac{1}{2401^s}\+\frac{2}{3721^s}\+\frac{4}{4225^s}\+\frac{2}{5329^s}
\+\frac{8}{6253^s}\+\frac{4}{8125^s}\+\frac{4}{8281^s}\+\frac{2}{9409^s}
\+\frac{8}{10309^s}\+\frac{2}{11881^s}
 +\ldots$\Strut\\
13 & $\frac{2}{729^s}\+\frac{6}{2809^s}\+\frac{6}{6241^s}\+\frac{6}{17161^s}
\+\frac{8}{19683^s}\+\frac{6}{24649^s}\+\frac{24}{38637^s}\+\frac{24}{57591^s}
\+\frac{24}{75843^s}\+\frac{24}{95499^s}\+\frac{6}{97969^s}
\+\frac{24}{114453^s}\+\frac{72}{148877^s}
+ \ldots$\Strut\\
15 & $\frac{1}{256^s}\+\frac{4}{961^s}\+\frac{4}{3721^s}\+\frac{2}{4096^s}
\+\frac{8}{7936^s}\+\frac{2}{14641^s}\+\frac{8}{15376^s}\+\frac{8}{15616^s}
\+\frac{4}{22801^s}\+\frac{32}{29791^s}\+\frac{4}{30976^s}\+\frac{4}{32761^s}
\+\frac{8}{38656^s}
+ \ldots$\Strut\\
16 & $\frac{4}{289^s}\+\frac{2}{2401^s}\+\frac{32}{4913^s}\+\frac{1}{6561^s}
\+\frac{4}{9409^s}\+\frac{4}{12769^s}\+\frac{16}{14161^s}\+\frac{8}{23409^s}
\+\frac{32}{28033^s}\+\frac{32}{32657^s}\+\frac{4}{37249^s}
\+\frac{16}{40817^s}\+\frac{32}{55777^s}
+ \ldots$\Strut\\
17 & $\frac{8}{10609^s}\+\frac{8}{18769^s}\+\frac{8}{57121^s}
\+\frac{8}{94249^s}\+\frac{8}{167281^s}\+\frac{8}{196249^s}
\+\frac{8}{375769^s}\+\frac{8}{418609^s}\+\frac{8}{844561^s}
\+\frac{8}{908209^s}\+\frac{8}{1042441^s}
+\ldots$\Strut\\
19 & $\frac{9}{36481^s}\+\frac{9}{52441^s}\+\frac{3}{117649^s}
\+\frac{9}{175561^s}\+\frac{9}{208849^s}\+\frac{9}{326041^s}
\+\frac{9}{418609^s}\+\frac{9}{579121^s}\+\frac{9}{1216609^s}
\+\frac{9}{1481089^s}\+\frac{3}{1771561^s}
+ \ldots$\Strut\\
20 & $\frac{1}{25^s}\+\frac{2}{125^s}\+\frac{3}{625^s}\+\frac{8}{1025^s}
\+\frac{8}{1525^s}\+\frac{4}{1681^s}\+\frac{2}{2025^s}\+\frac{8}{2525^s}
\+\frac{4}{3025^s}\+\frac{4}{3125^s}\+\frac{4}{3721^s}\+\frac{8}{4525^s}
\+\frac{16}{5125^s}\+\frac{8}{6025^s}
+ \ldots$\Strut\\
21 & $\frac{1}{49^s}\+\frac{2}{343^s}\+\frac{6}{1849^s}\+\frac{12}{2107^s}
\+\frac{3}{2401^s}\+\frac{2}{3136^s}\+\frac{1}{4096^s}\+\frac{12}{6223^s}
\+\frac{6}{8281^s}\+\frac{12}{10339^s}\+\frac{12}{12943^s}
\+\frac{24}{14749^s}\+\frac{6}{16129^s}\+\frac{12}{16513^s}
+ \ldots$\Strut\\
24 & $\frac{1}{81^s}\+\frac{2}{625^s}\+\frac{2}{729^s}\+\frac{4}{2025^s}
\+\frac{2}{2401^s}\+\frac{4}{3969^s}\+\frac{4}{5329^s}\+\frac{4}{5625^s}
\+\frac{8}{5913^s}\+\frac{3}{6561^s}\+\frac{8}{7857^s}\+\frac{4}{9409^s}
\+\frac{4}{9801^s}\+\frac{4}{13689^s}
+ \ldots$\Strut\\
25 & $\frac{10}{10201^s}\+\frac{10}{22801^s}\+\frac{10}{63001^s}
\+\frac{10}{160801^s}\+\frac{10}{361201^s}\+\frac{10}{491401^s}
\+\frac{10}{564001^s}\+\frac{200}{1030301^s}\+\frac{10}{1104601^s}
\+\frac{10}{1324801^s}\+\frac{10}{1442401^s}
+ \ldots$\Strut\\
27 & $\frac{9}{11881^s}\+\frac{9}{26569^s}\+\frac{9}{73441^s}
\+\frac{9}{143641^s}\+\frac{9}{187489^s}\+\frac{9}{237169^s}
\+\frac{9}{292681^s}\+\frac{9}{573049^s}\+\frac{9}{657721^s}
\+\frac{9}{844561^s}\+\frac{162}{1295029^s}
+ \ldots$\Strut\\
28 & $\frac{1}{64^s}\+\frac{2}{512^s}\+\frac{6}{841^s}\+\frac{12}{1856^s}
\+\frac{3}{4096^s}\+\frac{12}{6728^s}\+\frac{12}{7232^s}\+\frac{6}{10816^s}
\+\frac{12}{12608^s}\+\frac{6}{12769^s}\+\frac{24}{14848^s}
\+\frac{12}{17984^s}\+\frac{12}{21568^s}\+\frac{72}{24389^s}
+ \ldots$\Strut\\
32 & $\frac{8}{9409^s}\+\frac{8}{37249^s}\+\frac{8}{66049^s}
\+\frac{4}{83521^s}\+\frac{8}{124609^s}\+\frac{8}{201601^s}
\+\frac{8}{332929^s}\+\frac{8}{410881^s}\+\frac{8}{452929^s}
\+\frac{8}{591361^s}\+\frac{8}{863041^s}\+\frac{128}{912673^s}
+ \ldots$\Strut\\
33 & $\frac{10}{4489^s}\+\frac{10}{39601^s}\+\frac{1}{59049^s}
\+\frac{10}{109561^s}\+\frac{10}{157609^s}\+\frac{10}{214369^s}
\+\frac{5}{279841^s}\+\frac{200}{300763^s}\+\frac{10}{436921^s}
\+\frac{10}{528529^s}\+\frac{10}{737881^s}\+\frac{200}{893311^s}
+ \ldots$\Strut\\
35 & $\frac{12}{5041^s}\+\frac{12}{44521^s}\+\frac{12}{78961^s}
\+\frac{12}{177241^s}\+\frac{12}{241081^s}\+\frac{288}{357911^s}
\+\frac{12}{398161^s}\+\frac{12}{491401^s}\+\frac{6}{707281^s}
\+\frac{12}{829921^s}\+\frac{288}{1063651^s}
+ \ldots$\Strut\\
36 & $\frac{6}{1369^s}\+\frac{6}{5329^s}\+\frac{6}{11881^s}
\+\frac{6}{32761^s}\+\frac{72}{50653^s}\+\frac{3}{83521^s}
\+\frac{72}{99937^s}\+\frac{3}{130321^s}\+\frac{72}{149221^s}
\+\frac{6}{157609^s}\+\frac{6}{187489^s}\+\frac{72}{197173^s}
+ \ldots$\Strut\\
40 & $\frac{1}{625^s}\+\frac{8}{1681^s}\+\frac{2}{6561^s}
\+\frac{4}{14641^s}\+\frac{2}{15625^s}\+\frac{16}{25625^s}
\+\frac{16}{42025^s}\+\frac{4}{50625^s}\+\frac{8}{58081^s}
\+\frac{128}{68921^s}\+\frac{8}{75625^s}\+\frac{8}{78961^s}
+ \ldots$\Strut\\
44 & $\frac{10}{7921^s}\+\frac{10}{124609^s}\+\frac{10}{157609^s}
\+\frac{5}{279841^s}\+\frac{10}{380689^s}\+\frac{10}{436921^s}
\+\frac{200}{704969^s}\+\frac{10}{776161^s}\+\frac{10}{1026169^s}
\+\frac{10}{1630729^s}\+\frac{10}{1745041^s}
+ \ldots$\Strut\\
45 & $\frac{12}{32761^s}\+\frac{12}{73441^s}\+\frac{6}{130321^s}
\+\frac{12}{292681^s}\+\frac{12}{398161^s}\+\frac{12}{657721^s}
\+\frac{12}{982081^s}\+\frac{12}{1371241^s}
\+\frac{12}{2343961^s}\+\frac{12}{2627641^s}\+\frac{12}{3243601^s}
+  \ldots$\Strut\\
48 & $\frac{4}{2401^s}\+\frac{1}{6561^s}\+\frac{8}{9409^s}\+\frac{8}{37249^s}
\+\frac{8}{58081^s}\+\frac{4}{83521^s}\+\frac{8}{113569^s}\+\frac{32}{117649^s}
\+\frac{8}{187489^s}\+\frac{8}{194481^s}\+\frac{64}{232897^s}
\+\frac{4}{279841^s}
+ \ldots$\Strut\\
60 & $\frac{1}{256^s}\+\frac{1}{625^s}\+\frac{8}{3721^s}\+\frac{2}{4096^s}
\+\frac{2}{6400^s}\+\frac{1}{6561^s}\+\frac{2}{10000^s}\+\frac{4}{14641^s}
\+\frac{16}{15616^s}\+\frac{2}{15625^s}\+\frac{2}{20736^s}\+\frac{8}{30976^s}
\+\frac{8}{32761^s}
+ \ldots$\Strut\\
84 & $\frac{1}{2401^s}\+\frac{1}{4096^s}\+\frac{6}{28561^s}\+\frac{12}{113569^s}
\+\frac{2}{117649^s}\+\frac{2}{153664^s}\+\frac{12}{177241^s}\+\frac{2}{200704^s}
\+\frac{2}{262144^s}\+\frac{12}{405769^s}\+\frac{12}{452929^s}\+\frac{12}{573049^s}
+ \ldots$\Strut\\[2mm] 
\hline
\end{tabular*}
\renewcommand{\arraystretch}{1}
\end{table*}

Comparing $\varPsi_n$ with $\varPhi_n$, one notices that the two Dirichlet
series have precisely the same terms non-vanishing, i.e., $b_{n}(k)\neq 0$
if and only if $c^{}_{n}(k)\neq 0$. This implies

\begin{coro} {\rm \bf (CN1)}.
   The spectrum of possible coincidence indices remains unchanged by the
   addition of multiple coincidences. 
   In particular, the total spectrum\/ $\Sigma_{\OO_{n}}$ is the semigroup
   generated by\/ $1$ and the basic indices given in Tables\/ $\ref{tab:basic}$ 
   and\/ $\ref{tab:rami}$.
\qed
\end{coro} 

Moreover, as in Section 2, all multiple CSLs or CSMs can actually be
obtained from double intersections.  Consequently, $\varPsi_{n}(s)$ is,
once again, also the generating function for this situation.

Clearly, there is an ordering, $\varPsi_{n}(s)\succcurlyeq \varPhi_{n}(s)$, 
in the sense that $b_{n}(k)\ge c^{}_{n}(k)$ for all $k\ge 1$. In
Table~\ref{diritab2}, we list the first terms of the difference,
$\varPsi_{n}(s)-\varPhi_{n}(s)$. In fact, one quickly checks that
\begin{equation} \label{mono}
  \zeta^{}_{\KK_n}(s) \; \succcurlyeq \; \varPsi_{n}(s)
  \; \succcurlyeq \; \varPhi_{n}(s)
\end{equation}
in this sense, because every CSM of $\OO_n$ is a principal ideal of 
this ring, but not necessarily vice versa, so that
\[
    a_{n} (k) \; \ge \; b_{n} (k) \; \ge \; c^{}_{n} (k)
\]
for all $k\in\NN$. The summatory functions for 
the coefficients of $\varPhi_n (s)$ and of $\zeta^{}_{\KK_n}(s)$ both
show linear growth, see Corollary \ref{coro1} and \cite[Thm.~4]{BG2}.
Consequently, by monotonicity, the same behaviour shows up for
the coefficients of $\varPsi_n (s)$. 

To gain a better understanding of $\varPsi_{n}(s)$, and to prepare for
our later analysis of asymptotic properties, we observe the fundamental
relation
\begin{equation} \label{step1}
    \varPsi_{n}(s) \; = \; 
    \varPhi_{n}(s) \cdot \big\{ \varPsi_{n} (2s)\big\}^{1/2},
\end{equation}
which is certainly valid for all $s\in\CC$ with
$\mathrm{Re}(s)>1$. This relation can directly be derived from the
Euler product expansions of $\varPsi_{n}$ and $\varPhi_{n}$, compare
\cite[Sec.~2]{Sh} for a similar situation and \cite{MC} for a
historical perspective.

Moreover, iterating Eq.~\eqref{step1} leads to
\begin{prop}\label{prop2}
  The two Dirichlet series generating functions\/ $\varPhi_{n}(s)$ and\/
  $\varPsi_{n}(s)$ of Theorems\/ $\ref{thm1}$ and\/ $\ref{thm2}$ are
  related by
\[
  \varPsi_{n}(s) \; = \; \big(\varPsi_{n}(2^{L+1} s)\big)^{1/2^{L+1}}
  \cdot \prod_{\ell=0}^{L}
  \big( \varPhi_{n}(2^{\ell}s)\big)^{1/2^{\ell}},
\]
  for any integer\/ $L\ge 0$. Moreover, one also has
\[
  \varPsi_{n}(s) \; = \; \prod_{\ell=0}^{\infty}
  \big( \varPhi_{n}(2^{\ell}s)\big)^{1/2^{\ell}}
   = \; \varPhi_{n}(s) \cdot  \prod_{\ell=1}^{\infty}
  \big( \varPhi_{n}(2^{\ell}s)\big)^{1/2^{\ell}},
\]
  which is absolutely convergent on\/ $\{ \mathrm{Re}(s) > 1\}$.
\end{prop}

\noindent {\sc Proof}.  
The first claim follows by simple induction (in $L$)
from Eq.~\eqref{step1}.

If $\mathrm{Re}(s) > 1$, it is clear that
$\mathrm{Re} (2^L s) \longrightarrow\infty$ as $L\to\infty$.
In this limit, by well-known properties of convergent
Dirichlet series, one has
\[
   \varPsi_n (2^L s) \longrightarrow b_n (1) = 1 \quad\mbox{and}\quad
   \varPhi_n (2^L s) \longrightarrow c_n (1) = 1 \, .
\]
If we write $s=\sigma + i t$, these limits are uniform in $t\in\RR$,
see \cite[Thm.~11.2]{A}. One obvious consequence is
\[
   \lim_{L\to\infty}
   \left( \varPsi_n (2^L s) \right)^{1/2^L} \; = \; 1 \, ,
\]
but the above observation also gives the absolute convergence
of the infinite product involved, because its logarithm,
\[
   \sum_{\ell=0}^{\infty} \frac{1}{2^\ell}
   \log\big({\varPhi_n (2^\ell s)}\big),
\]
converges absolutely on $\{\mathrm{Re}(s)>1\}$.
\qed\smallskip

Observe that, by analytic continuation, the product
\[
    \prod_{\ell=1}^{\infty}
    \big( \varPhi_{n}(2^{\ell}s)\big)^{1/2^{\ell}},
\]
starting with $\varPhi_n (2s)^{1/2}$,
defines an analytic function on $\{ \mathrm{Re}(s) > 1/2 \}$ without
any zeros in this half-plane (the absence of zeros can actually be extended 
to its closure, $\{ \mathrm{Re}(s) \ge 1/2 \}$, by known properties of
$\zeta(s)$ and the $L$-series involved). This shows that $\varPhi_{n}(s)$
and $\varPsi_{n}(s)$ share the positions of their zeros and poles, when 
defined by analytic continuation on this half-plane. In particular, they 
both have a single simple pole in $\{ \mathrm{Re}(s) > 1/2 \}$, at $s=1$,
and they share the position (if any) and orders of zeros there. This
remains true on the critical line $\{ \mathrm{Re}(s) = 1/2 \}$, where
the only zeros of $\varPhi_{n}(s)$ (and hence of $\varPsi_{n}(s)$) 
in the entire set $\{ \mathrm{Re}(s) \ge 1/2 \}$ are expected on the
basis of the Riemann hypothesis for $\zeta(s)$ and its generalisation 
to $L$-series.

Without any reference to this (still unproved) hypothesis, one has,
by standard arguments on the basis of compact convergence,
\begin{coro} \label{quick-conv}
  On the half-plane\/ $\{ \mathrm{Re}(s) > 1/2 \}$, by unique 
  analytic continuation, one has the representation
\[
   \frac{\varPsi_{n}}{\varPhi_{n}}(s) \; = \; 
   \prod_{\ell=1}^{\infty} \big(\varPhi_{n} (2^{\ell} s)
                  \big)^{1/2^\ell}
\]
  which is an analytic and zero-free function in this half-plane.

  In particular, $\varPsi_{n}/\varPhi_{n}$ has a well-defined
  positive value at\/ $s=1$.    \qed
\end{coro}

At this point, we can state the asymptotic result for the
summatory function of the coefficients $b_n (k)$.
\begin{coro}
{\rm \bf (CN1)}. 
The number of  multiple CSMs of index\/ $\le x$ 
shows the asymptotic behaviour
\[
  \sum_{k\le x} b_{n}(k) \; \sim \; x\cdot 
  \big(\mbox{\rm res}_{s=1}\,\varPsi_{n}(s)\big)
  \; = \;  x \cdot \beta_n ,
\]
with the growth constant
\[
  \beta_n \; = \;  q_n \cdot \big(\mbox{\rm res}_{s=1}\,\varPhi_{n}(s)\big)
          \; = \;  q_n\cdot \gamma_n .
\]
Here, the constant $q_n$ has the monotonically increasing and
rapidly converging product representation
\[
   q_n \; := \; \lim_{s\to 1} \frac{\varPsi_n(s)}{\varPhi_n(s)}
    \; = \; \prod_{\ell=1}^{\infty} \big(\varPhi_n(2^\ell)\big)^{1/2^\ell}.
\]
\end{coro}
\noindent {\sc Proof}.
The monotonicity w.r.t.\ $\succcurlyeq$ stated in \eqref{mono}, together
with the analyticity properties of $\varPsi_n (s)$,  implies that 
$\varPsi_n (s)$ also has a simple pole at $s=1$, and no other 
singularity in $\{\mathrm{Re}(s)\ge 1\}$.
Delange's theorem then yields the claim on the asymptotic behaviour.

The second statement is almost immediate from our above discussion.
Clearly, $\varPhi_n (2^\ell)>1$ for all $\ell\in\NN$, hence also
$\varPhi_n (2^\ell)^{1/2^\ell} > 1$, which gives the
monotonicity, while convergence follows from Corollary~\ref{quick-conv} 
and Proposition~\ref{prop2}.
\qed\smallskip

In Table~\ref{tab:asym}, we give the numerical values of the residues
of the three types of generating functions we have encountered. Some
values for $\alpha_{n}$ and $\gamma_{n}$, including exact
expressions (except for $\gamma^{}_{7}$, see a previous footnote), 
are also contained in \cite{BG2} and \cite{PBR}, respectively.

\section{Outlook}

It is desirable to extend the above analysis to higher dimensions,
which is significantly more involved due to non-commutativity of
the $\mathrm{SOC}$-groups in these cases.
Nevertheless, quite a bit is known for simple
coincidences \cite{B,P1,P2}, and first steps are in sight for multiple
ones \cite{Z}.

Also, in the planar case, one would like to get rid of the
assumption ({\bf CN1}). This is already pretty tricky for simple
coincidences, see \cite{PBR} for an example ($n=23$), but some
further results seem possible.

Finally, there are many open questions concerning the
structure of the $\mathrm{SOC}$-groups, particularly for
dimensions $d\ge 3$. As they are a natural extension of
point symmetry groups, they certainly deserve further
attention.

\begin{acknowledgment}
It is a pleasure to thank Peter A.\ B.\ Pleasants and Peter Zeiner
for cooperation and helpful discussions, and P.\ Moree for useful
hints on the literature. This work was supported by the German
Research Council, within the Collaborative Research Centre 701.
\end{acknowledgment}


\begin{thebibliography}{99}\itemsep 0pt

\bibitem{A}
T.~M.~Apostol,
{\em An Introduction to Analytic Number Theory},
5th corr.\ printing (Springer, New York, 1998).

\bibitem{B}
M.~Baake,
Solution of the coincidence problem in dimensions \mbox{$d\le 4$},
in: {\em The Mathematics of Longe-Range Aperiodic Order},
ed.\ R.\ V.\ Moody, NATO ASI Series C 489
(Kluwer, Dordrecht, 1997) pp.\ 9--44;
rev.\ version, {\tt math.MG/0605222}.

\bibitem{B2}
M.~Baake,
A guide to mathematical quasicrystals, in:
{\em Quasicrystals. 
An Introduction to Structure, Physical Properties and Applications},
eds.\ J.-B.\ Suck, M.\ Schreiber and P.\ H\"{a}ussler 
(Springer, Berlin, 2002) pp.\ 17--48;
{\tt math-ph/9901014}.

\bibitem{BG}
M.~Baake and U.~Grimm,
Combinatorial problems of (quasi)\-crystallography,
in: {\em Quasicrystals:\ Structure and Physical
Properties}, ed.\ H.-R.\ Trebin
(Wiley-VCH, Weinheim, 2003) pp.\ 160--171;
{\tt math-ph/0212015}.

\bibitem{DCG}
M.~Baake and U.~Grimm,
A note on shelling, 
{\em Discr.\ Comput.\ Geom.}\/ {\bf 30} (2003) 573--589; 
{\tt math.MG/0203025}.

\bibitem{BG2}
M.~Baake and U.~Grimm,
Bravais colourings of planar modules with $N$-fold symmetry,
{\em Z.\ Kristallographie}\/ {\bf 219} (2004) 72--80; 
{\tt math.CO/0301021}.

\bibitem{BP}
M.~Baake and P.~A.~B.~Pleasants,
Algebraic solution of the coincidence problem in two
and three dimensions,
{\em Z.\ Naturf.}\/ {\bf 50a} (1995) 711--717.

\bibitem{Berndt}
B.~C.~Berndt, On the Hurwitz zeta-function,
{\em Rocky Mountain J.\ Math.}\/ {\bf 2} (1972) 151--157.

\bibitem{Bollmann}
W.~Bollmann,
{\em Crystal Defects and Crystalline Interfaces}\/
(Springer, Berlin, 1970).

\bibitem{BS} 
Z.~I.~Borevich and I.~R.~Shafarevich,
{\em Number Theory}\/ (Academic Press, New York, 1966).

\bibitem{Gertsman}
V.~Y.~Gertsman, 
Geometrical theory of triple junctions of CSL boundaries,
{\em Acta Cryst.}\/ {\bf A 57} (2001) 369--377.

\bibitem{GP}
U.~Grimm,\newline
{\tt http://mcs.open.ac.uk/ugg2/coincidences/}.

\bibitem{G}
H.~Grimmer,
Coincidence orientations of grains in rhombohedral materials,
{\em Acta Cryst.}\/ {\bf A 45} (1989) 505--523.

\bibitem{HW}
G.~H.~Hardy and E.~M.~Wright,
{\em An Introduction to the Theory of Numbers}, 5th ed.
(Clarendon Press, Oxford 1979).
 
\bibitem{MC}
P.~Moree and J.~Cazaran,
On a claim of Ramanujan in his first letter to Hardy,
{\em Exp.\ Math.}\/ {\bf 17} (1999) 289--312.

\bibitem{P}
P.~A.~B.~Pleasants, 
Designer quasicrystals:\ Cut-and-project sets
with pre-assigned properties,
in:\  {\em Directions in Mathematical Quasicrystals},
eds.\ M.\ Baake and R.\ V.\ Moody,
CRM Monograph Series, vol.~13 
(AMS, Providence, RI, 2000) pp.~95--141.

\bibitem{PBR}
P.~A.~B.~Pleasants, M.~Baake and J.~Roth,
Planar coincidences for $N$-fold symmetry,
{\em J.\ Math.\ Phys.}\/ {\bf 37} (1996) 1029--1058; 
{\tt math.MG/0511147}.

\bibitem{Sh}
D.~Shanks,
The second-order term in the asymptotic expansion of $B(x)$,
{\em Math.\ Comp.}\/ {\bf 18} (1964) 78--86.

\bibitem{Sloane}
N.~J.~E.~Sloane and B.~Beferull-Lozano,
Quantizing using lattice intersections,
in:\ {\em Discrete and Computational Geometry},
eds.\ B.\ Aronov, S.\ Basu, J.\ Pach and M.\ Sharir
(Springer, Berlin, 2003) pp.~799--824;
{\tt math.CO/0207147}.

\bibitem{Steurer}
W.~Steurer,
Twenty years of structure research on quasicrystals.
Part I. Pentagonal, octagonal, decagonal and
dodecagonal quasicrystals,
{\em Z.\ Kristallographie} {\bf 219} (2004) 391--446.

\bibitem{T}
G.~Tenenbaum, 
{\em Introduction to Analytic and Probabilistic Number Theory}\/ 
(Cambridge University Press, Cambridge, 1995).

\bibitem{W}
L.~C.~Washington, 
{\em Introduction to Cyclotomic Fields}, 2nd ed.\
(Springer, New York, 1997).

\bibitem{WW}
E.~T.~Whittaker and G.~N.~Watson,
{\em A Course in Modern Analysis}, 4th ed.\ (reprinted)
(Cambridge University Press, Cambridge, 1965).

\bibitem{P1}
P.~Zeiner,
Symmetries of coincidence site lattices of cubic lattices,
{\em Z.\ Kristallographie} {\bf 220} (2005) 915--925.

\bibitem{P2}
P.~Zeiner,
Coincidences of hypercubic lattices in 4 dimensions
{\em Z.\ Kristallographie} {\bf 221} (2006) 105--114.

\bibitem{Z}
P.~Zeiner, 
Multiple CSLs for the body centered cubic lattice,
preprint (2006).

\end{thebibliography}
\end{document}